\documentclass{amsproc}

\theoremstyle{definition}

\theoremstyle{remark}

\numberwithin{equation}{section}

\def\E{\ifmmode{\mathbb E}\else{$\mathbb E$}\fi} 
\def\N{\ifmmode{\mathbb N}\else{$\mathbb N$}\fi} 
\def\R{\ifmmode{\mathbb R}\else{$\mathbb R$}\fi} 
\def\Q{\ifmmode{\mathbb Q}\else{$\mathbb Q$}\fi} 
\def\C{\ifmmode{\mathbb C}\else{$\mathbb C$}\fi} 
\def\H{\ifmmode{\mathbb H}\else{$\mathbb H$}\fi} 
\def\Z{\ifmmode{\mathbb Z}\else{$\mathbb Z$}\fi} 
\def\P{\ifmmode{\mathbb P}\else{$\mathbb P$}\fi} 
\def\T{\ifmmode{\mathbb T}\else{$\mathbb T$}\fi} 
\def\SS{\ifmmode{\mathbb S}\else{$\mathbb S$}\fi} 
\def\DD{\ifmmode{\mathbb D}\else{$\mathbb D$}\fi} 

\newcommand{\e}{\varepsilon}

\newcommand{\del}{\partial}

\newcommand{\ben}{\begin{enumerate}}
\newcommand{\een}{\end{enumerate}}
\newcommand{\be}{\begin{equation}}
\newcommand{\ee}{\end{equation}}
\newcommand{\bea}{\begin{eqnarray}}
\newcommand{\eea}{\end{eqnarray}}
\newcommand{\beastar}{\begin{eqnarray*}}
\newcommand{\eeastar}{\end{eqnarray*}}
\newcommand{\bc}{\begin{center}}
\newcommand{\ec}{\end{center}}

\theoremstyle{theorem}
\newtheorem{thm}{Theorem}[section]
\newtheorem{cor}[thm]{Corollary}
\newtheorem{lem}[thm]{Lemma}
\newtheorem{prop}[thm]{Proposition}

\newtheorem{conj}[thm]{Conjecture}

\theoremstyle{definition}
\newtheorem{defn}[thm]{Definition}

\newtheorem{ques}[thm]{Question}

\newtheorem{exm}[thm]{Example}

\newtheorem*{thm*}{Theorem}

\numberwithin{equation}{section}

\def\R{{\mathbb R}}
\def\osc{{\hbox{\rm osc}}}

\def\E{{\mathbb E}}
\def\Z{{\mathbb Z}}
\def\C{{\mathbb C}}
\def\R{{\mathbb R}}
\def\P{{\mathbb P}}

\def\N{{\mathbb N}}

\def\11{{\mathbb I}}

\def\C{\mathbb{C}}
\def\Z{\mathbb{Z}}

\def\T{\mathbb{T}}

\def\Q{\mathbb{Q}}

\def\E{\ifmmode{\mathbb E}\else{$\mathbb E$}\fi} 
\def\N{\ifmmode{\mathbb N}\else{$\mathbb N$}\fi} 
\def\R{\ifmmode{\mathbb R}\else{$\mathbb R$}\fi} 
\def\Q{\ifmmode{\mathbb Q}\else{$\mathbb Q$}\fi} 
\def\C{\ifmmode{\mathbb C}\else{$\mathbb C$}\fi} 
\def\H{\ifmmode{\mathbb H}\else{$\mathbb H$}\fi} 
\def\Z{\ifmmode{\mathbb Z}\else{$\mathbb Z$}\fi} 
\def\P{\ifmmode{\mathbb P}\else{$\mathbb P$}\fi} 
\def\SS{\ifmmode{\mathbb S}\else{$\mathbb S$}\fi} 
\def\DD{\ifmmode{\mathbb D}\else{$\mathbb D$}\fi} 

\def\R{{\mathbb R}}
\def\osc{{\hbox{\rm osc}}}

\def\E{{\mathbb E}}
\def\Z{{\mathbb Z}}
\def\C{{\mathbb C}}
\def\R{{\mathbb R}}

\def\N{{\mathbb N}}

\def\leng{{\operatorname{leng}}}

\def\e{\varepsilon}

\def\CA{{\mathcal A}}

\def\CH{{\mathcal H}}

\def\CP{{\mathcal P}}
\def\CQ{{\mathcal Q}}

\def\CP{{\mathcal P}}

\def\CT{{\mathcal T}}
\def\CU{{\mathcal U}}

\def\darr#1{\raise1.5ex\hbox{$\leftrightarrow$}
\mkern-16.5mu #1}

\def\roughly#1{\raise.3ex\hbox{$#1$\kern-.75em
\lower1ex\hbox{$\sim$}}}

\def\opname#1{\mathop{\kern0pt{\rm #1}}\nolimits}

\begin{document}
\quad \vskip1.375truein

\def\mq{\mathfrak{q}}
\def\mp{\mathfrak{p}}
\def\mH{\mathfrak{H}}
\def\mh{\mathfrak{h}}
\def\ma{\mathfrak{a}}
\def\ms{\mathfrak{s}}
\def\mm{\mathfrak{m}}
\def\mn{\mathfrak{n}}

\def\Hoch{{\tt Hoch}}
\def\mt{\mathfrak{t}}
\def\ml{\mathfrak{l}}
\def\mT{\mathfrak{T}}
\def\mL{\mathfrak{L}}
\def\mg{\mathfrak{g}}
\def\md{\mathfrak{d}}


\title[Continuous Hamiltonian flows]{The group of Hamiltonian
homeomorphisms and continuous Hamiltonian flows}

\author{Yong-Geun Oh}
\address{Department of Mathematics, University of Wisconsin, Madison, WI
53706}
\email{oh@math.wisc.edu}

\thanks{Partially supported by the NSF grant \# DMS 0503954 and a grant of the
2000 Korean Young Scientist Prize}


%
%

\begin{abstract}
In this paper, we study the dynamical aspects of the
\emph{Hamiltonian homeomorphism group} $Hameo(M,\omega)$ which was
introduced by M\"uller and the author. We introduce the notion of
autonomous continuous Hamiltonian flows and extend the well-known
conservation of energy to such flows. The definitions of the Hofer
length and of the spectral invariants $\rho_a$ are extended to
continuous Hamiltonian paths, and the Hofer norm and the spectral
norm $\gamma:Ham(M,\omega) \to \R_+$ are generalized to the
corresponding intrinsic norms on $Hameo(M,\omega)$ respectively.
Using these extensions, we also extend the construction of
Entov-Polterovich's Calabi quasi-morphism on $S^2$ to the space of
continuous Hamiltonian paths. We also discuss a conjecture
concerning extendability of Entov-Polterovich's quasi-morphism and
its relation to the extendability of Calabi homomorphism on the disc
to $Hameo(D^2,\del D^2)$, and their implication towards the
simpleness question on the area preserving homeomorphism groups of
the disc $D^2$ and of the sphere $S^2$.
\end{abstract}

\maketitle

\section{Introduction}
\label{sec:intro}

\subsection{Topological Hamiltonian flows}
\label{subsec:top-flows}

A time-dependent Hamilton's equation on a symplectic manifold
$(M,\omega)$ is the first order ordinary differential equation
$$
\dot x = X_H(t,x)
$$
where the time-dependent vector field $X_H$ associated to a function
$H: \R \times M \to \R$ is given by the defining equation
\be\label{eq:defining} dH_t = X_{H_t}\rfloor \omega. \ee Therefore
if we consider functions $H$ \emph{that are $C^{1,1}$ so that one
can apply the existence and uniqueness theorem of solutions of the
above Hamilton's equation}, the flow $t \mapsto \phi_H^t$, an
isotopy of diffeomorphisms, is uniquely determined by the
Hamiltonian $H$. We will always assume
\begin{enumerate}
\item the Hamiltonians are normalized by $\int_M H_t \, d\mu = 0$
for the Liouville measure $d\mu$ of $(M,\omega)$ if $M$ is closed,
\item and they are compactly supported in $\mbox{Int}M$ if $M$ is
open.
\end{enumerate}
We call such Hamiltonian functions {\it normalized}. For the
convenience of exposition, we will focus on the closed case unless
otherwise said. All the discussions in this paper equally apply for
the open case too.

We denote by $C_m^\infty(M)$ the set of normalized smooth functions
on $M$ and by $\CP(C_m^\infty(M))= C^\infty_m([0,1] \times M)$ the
set of time-dependent normalized Hamiltonian functions. We will also
denote the Hamiltonian isotopy generated by $H$ by
$$
\phi_H: t\mapsto \phi_H^t.
$$
Conversely if a smooth isotopy $\lambda$ of Hamiltonian
diffeomorphisms is given, we can obtain the corresponding normalized
Hamiltonian $H$ by differentiating the isotopy and then solving
(\ref{eq:defining}). Therefore \emph{in the smooth category} this
correspondence is bijective.

On the other hand, due to the fact that this correspondence involves
differentiating the function and solving Hamilton's equation, the
correspondence gets murkier as the regularity of the Hamiltonian is
weaker than $C^{1,1}$ because of solvability question of Hamilton's
equation.

In \cite{oh:hameo1}, the author and M\"uller studied this relation
and introduced the notion of \emph{hamiltonian limits} of smooth
Hamiltonian flows and proposed the notion of \emph{continuous
Hamiltonian flow} as the hamiltonian limits thereof. Then we
introduced the $C^0$-concept of Hamiltonian diffeomorphisms, called
\emph{Hamiltonian homeomorphisms}, which forms a normal subgroup of
the group of \emph{symplectic homeomorphisms}: Motivated by
Eliashberg's $C^0$-symplectic rigidity theorem \cite{eliash}, we
defined in \cite{oh:hameo1} the group of \emph{symplectic
homeomorphisms} as follows. We give the compact-open topology on
$Homeo(M)$, which is equivalent to the metric topology induced by
the metric
$$
\overline d(\phi,\psi) = \max\{d_{C^0}(\phi,\psi), d_{C^0}
(\phi^{-1},\psi^{-1})\}
$$
on a compact manifold $M$.

\begin{defn}[Symplectic homeomorphism group]
Define $Sympeo(M,\omega)$ to be
$$
Sympeo(M,\omega):= \overline{Symp(M,\omega)}
$$
the $C^0$-closure of $Symp(M,\omega)$ in $Homeo(M)$ and call
$Sympeo(M,\omega)$ the {\it symplectic homeomorphism group}.
\end{defn}

We now recall the formal definition of continuous Hamiltonian flow
introduced in \cite{oh:hameo1}. Hofer's $L^{(1,\infty)}$ norm of
Hamiltonian diffeomorphisms is defined by
$$
\|\phi\| = \inf_{H\mapsto \phi} \|H\|
$$
where $H \mapsto \phi$ means that $\phi= \phi_H^1$ is the time-one
map of Hamilton's equation
$$
\dot x = X_H(t,x)
$$
and the norm $\|H\|$ is  defined by \be\label{eq:intosc} \|H\| =
\int _0^1 \mbox{\rm osc }H_t \, dt = \int_0^1(\max_x H_t -\min_x
H_t)\, dt. \ee For two given continuous paths $\lambda, \, \mu:[a,b]
\to Homeo(M)$, we define their distance by \be\label{eq:bard}
\overline d(\lambda,\mu) = \max_{t\in [a,b]}\overline
d(\lambda(t),\mu(t)). \ee Following \cite{oh:hameo1}, we denote by
$$
\CP^{ham}(Symp(M,\omega),id)
$$
the set of smooth Hamiltonian paths $\lambda: [0,1] \to
Symp(M,\omega)$ with $\lambda(0) = id$, and equip it with the
\emph{Hamiltonian topology} \cite{oh:hameo1}.

\begin{defn}[$C^0$-Hamiltonian topology]
Let $(M,\omega)$ be a closed symplectic manifold.
\begin{enumerate}
\item We define the \emph{$C^0$-Hamiltonian topology} of
the set $\CP^{ham}(Symp(M,\omega),id)$ of Hamiltonian paths  by the
one generated by the collection of subsets \bea
\begin{split}
& \CU(\phi_H,\e_1,\e_2):= \\
& \Big\{\phi_{H'} \in \CP^{ham}(Symp(M,\omega),id) \Big| \|\overline
H \# H'\| < \e_1, \, \overline d(\phi_H, \phi_{H'}) < \e_2 \Big\}
\end{split}
\eea of $\CP^{ham}(Symp(M,\omega),id)$ for $\e_1, \, \e_2> 0 $ and
$\phi_H \in \CP^{ham}(Symp(M,\omega),id)$. We denote the resulting
topological space by $\CP^{ham}_s(Symp(M,\omega),id)$.
\item
We define the \emph{$C^0$-Hamiltonian topology} of $Ham(M,\omega)$
to be the strongest topology such that the evaluation map
\be\label{ev1} ev_1: \CP^{ham}(Symp(M,\omega),id) \to Ham(M) \ee is
continuous. We denote the resulting topological space by $\CH
am(M,\omega)$.
\end{enumerate}
We will call continuous maps with respect to the Hamiltonian
topology \emph{Hamiltonian continuous}.
\end{defn}
The $C^0$-Hamiltonian topology of $\CP^{ham}(Symp(M,\omega),id)$ is
equivalent to the metric topology induced by the metric
$$
d_{ham}(\lambda,\mu): = \overline d(\lambda,\mu) +
\operatorname{leng}(\lambda^{-1}\mu)
$$
where $\overline d$ is the $C^0$ metric on $\CP(Homeo(M),id)$. (See
Proposition 3.10 \cite{oh:hameo1}.)

\begin{defn}[Topological Hamiltonian flow]\label{topflowdefn} A continuous map
$\lambda: \R \to Homeo(M)$ is called a topological Hamiltonian flow
if there exists a sequence of smooth Hamiltonians $H_i: \R \times M
\to \R$ satisfying the following:
\begin{enumerate}
\item $\phi_{H_i} \to \lambda$ locally uniformly on $\R \times M$.
\item the sequence $H_i$ is Cauchy in the $L^{(1,\infty)}$-topology and so
has a limit $H_\infty$ lying in $L^{(1,\infty)}$.
\end{enumerate}
We call a continuous path $\lambda:[a,b] \to Homeo(M)$ a {\it
topological Hamiltonian path} if it satisfies the same conditions
with $\R$ replaced by $[a,b]$, and the limit
$L^{(1,\infty)}$-function $H_\infty$ a \emph{topological
Hamiltonian}. In any of these cases, we say that the pair
$(\lambda,H_\infty)$ is the \emph{hamiltonian limit} of
$(\phi_{H_i},H_i)$, and write
$$
\operatorname{hlim}_{i \to \infty} (\phi_{H_i},H_i) \to (\lambda,
H_\infty)
$$
or sometimes even $\operatorname{hlim}_{i \to \infty}
(\phi_{H_i},H_i)=\lambda$.
\end{defn}

We denote by $\CP^{ham}_{[a,b]}(Sympeo(M,\omega),id)$ the set of
topological Hamiltonian paths defined on $[a,b]$. When $[a,b]=[0,1]$
or when we do not specify the domain of $\lambda$, we often just
write $\CP^{ham}(Sympeo(M,\omega),id)$ for the corresponding set of
topological Hamiltonian paths.

\begin{defn}[Hamiltonian homeomorphism group] We define
$$
Hameo(M,\omega) = ev_1(\CP_{[0,1]}^{ham}(Sympeo(M,\omega),id))
$$
and call any element therein a \emph{Hamiltonian homeomorphisms}
\end{defn}
One basic theorem proved in \cite{oh:hameo1} is that
$Hameo(M,\omega)$ forms a path-connected normal subgroup of
$Sympeo(M,\omega)$.

\subsection{Continuous Hamiltonian flows: statement of main results}
\label{subsec:cont-flows}

All the above discussion can be carried out using the stronger
version, or the $L^\infty$-version of Hofer's norm
$$
\|H\|_\infty: = \max_{t \in [0,1]}(\osc H_t)
$$
and define the set $\CP_\infty^{ham}(Sympeo(M,\omega),id)$ by
replacing $\| \cdot \|$ by $\|\cdot \|_{\infty}$. We call any
element thereof a {\it continuous Hamiltonian path}.

In \cite{oh:hameo1}, we defined the set of \emph{$C^0$-Hamiltonian
homeomorphisms} by \be\label{eq:hameo} Hameo_\infty(M,\omega) = \{h
\in Homeo(M) \mid h = \overline{ev}_1(\lambda), \, \lambda\in
\CP_\infty^{ham}(Sympeo(M,\omega),id)\}. \ee The following theorem
was proved by M\"uller \cite{mueller}
\begin{thm}[M\"uller]\label{muellerthm} We have
$$
Hameo(M,\omega) = Hameo_{\infty}(M,\omega).
$$
\end{thm}

In view of this theorem, we will drop $\infty$ from
$Hameo_\infty(M,\omega)$ from now on.

The following uniqueness theorem for continuous Hamiltonian flows,
which was asked in \cite{oh:hameo1} (in the
$L^{(1,\infty)}$-context), was proved by Viterbo for the closed
manifolds and by the author for the compactly supported case on open
manifolds

\begin{thm}[Viterbo \cite{viterbo2}, Oh \cite{oh:locality}]\label{thm:unique}
Suppose that $H_i$ is a sequence of smooth normalized Hamiltonians
that uniformly converges in $[0,1] \times M$. Then if $\phi_{H_i}$
converges uniformly to the constant path $id$, then we have $\lim_i
H_i \equiv 0$.
\end{thm}

We refer to \cite{oh:locality} for the locality result of
$C^0$-Hamiltonians. The following uniqueness result in the context
of $L^{(1,\infty)}$ Hamiltonians is still open.

\begin{ques} Does the uniqueness hold for $L^{(1,\infty)}$-Hamiltonians?
\end{ques}

Because of this lack of uniqueness result in the more natural
$L^{(1,\infty)}$-context, we will restrict our discussions from now
on in this paper mostly to the context of continuous Hamiltonian
flows. (However all the results would hold true and could be proven
in the same way for the topological Hamiltonian flows if this
uniqueness theorem should be proved in the $L^{(1,\infty)}$
context.)

Using the uniqueness result, we first establish the following
one-one correspondence. This extends the well-known correspondence
in the smooth category to this continuous Hamiltonian category. See
the next sections for more precise statements and some discussion on
this correspondence in perspective.

\begin{thm} We have a canonical one-one correspondence
\be\label{eq:correspond} \CH^0_m([0,1] \times M, \R)
\longleftrightarrow \CP^{ham}_\infty(Sympeo(M,\omega),id) \ee where
$\CH^0_m([0,1] \times M, \R)$ is the set of (normalized) continuous
Hamiltonians and $\CP^{ham}_\infty(Sympeo(M,\omega),id)$ is the set
of continuous Hamiltonian paths. And the following diagram commutes
: \be
\begin{matrix} C^{\infty}_m([0,1] \times M, \R) &\longleftrightarrow
& \CP^{ham}(Symp(M,\omega),id) \\
\downarrow & {} & \downarrow \\
\CH^0_m([0,1] \times M, \R) &\longleftrightarrow &
\CP^{ham}_\infty(Sympeo(M,\omega),id)
\end{matrix}
\ee where the vertical maps are canonical inclusion maps.
\end{thm}

We like to point out that the set $\CH^0_m([0,1] \times M, \R)$
contains all $C^{1,1}$ functions and
$\CP_\infty^{ham}(Sympeo(M,\omega),id)$ contains all the Hamiltonian
paths generated by $C^{1,1}$ functions.

The correspondence (\ref{eq:correspond}) can be interpreted as the
criterion for a $C^0$-Hamiltonian $H$ to allow a weak solution of
Hamilton's equation $\dot x = X_{H}(t,x)$. It would be interesting
to make this statement more precise in the point of view of the
generalized or distribution solutions of ordinary differential
equations.

As we will illustrate by several theorems concerning the general
properties of continuous Hamiltonian flows, this one-one
correspondence will be a crucial ingredient to prove those theorems.
We refer to Theorem \ref{one-one}, Theorem \ref{autonom},
Proposition \ref{hoferlength}, Theorem \ref{spectral} and Theorem
\ref{barmupath}, for example.

It is interesting to study this correspondence between continuous
Hamiltonian paths and their Hamiltonians in the point of view of Lie
group theory, which is a subject of future study. However, as a step
towards this goal, we prove the following theorem using the
uniqueness theorem. We refer to Theorem \ref{autonom} for a more
precise statement.

\begin{thm}
Suppose that $\lambda \in \CP_\infty^{ham}(Sympeo(M,\omega),id)$ is
a one-parameter subgroup, i.e., a path satisfying
$$
\lambda(t+s) = \lambda(t)\lambda(s)
$$
whenever $t,\, s$ and $t+s$ lie in the domain of $\lambda$. Then its
Hamiltonian $H$ is continuous and time-independent. And the converse
also holds.
\end{thm}

Similarly via the uniqueness theorem, we also prove in Theorem
\ref{conservation} that the autonomous continuous Hamiltonian flow
preserves its Hamiltonian.

\begin{thm} Let $H$ be an autonomous continuous Hamiltonian and
$\phi_H$ be its flow. Then we have $H\circ \phi_H^s = H$ for all
$s$.
\end{thm}
This generalizes the well-known conservation law of smooth
autonomous Hamiltonian mechanics. It leads us to the following
natural $C^0$-generalization of the existence question of a
time-periodic closed orbit for autonomous Hamiltonian system

\begin{ques} Does there always exist a periodic orbit of
continuous Hamiltonian flows on the hypersurface of a `generic'
level of a convex autonomous continuous Hamiltonian?
\end{ques}

The uniqueness theorem is also essential to extend the definitions
of the Hofer length and the spectral invariants constructed in
\cite{oh:alan} to continuous Hamiltonian paths: All the smooth
constructions concerning the Hofer length and the spectral
invariants use the Hamiltonian functions in their constructions, but
not directly their associated Hamiltonian paths. They are
interpreted as the invariants of the latter only via the one-one
correspondence between the Hamiltonian flows and the Hamiltonian
functions. Therefore extending these constructions to the
topological category attached to the continuous Hamiltonian
\emph{paths}, not to the \emph{functions}, requires this uniqueness
theorem. We refer readers to sections \ref{sec:hoferlength} and
\ref{sec:spectral} for the study of these extensions. In section
\ref{sec:quasimorphism}, as an application of this generalization of
spectral invariants, we extend Entov-Polterovich's study
\cite{entov-pol} of Calabi quasi-morphisms to the space
$\CP^{ham}(Sympeo(S^2),id)$ of continuous Hamiltonian paths on $S^2$
and state a conjecture (Conjecture \ref{barmu}). Finally in section
\ref{sec:wild}, we discuss an implication of this conjecture to the
simpleness question of area preserving homeomorphism group of the
disc and of the sphere.

We thank S. M\"uller and A. Fathi for many helpful comments and
discussions. We are also very grateful to the anonymous referee for
providing useful comments and suggestions to improve the
presentation and English of the paper and for asking us a question
the answer of which we provide by adding Theorem
\ref{thm:5.7imply5.8} and \ref{thm:simpleonS2} to this version of
the paper.

\bigskip

\noindent{\bf Notation:} Since our main concern in this paper is in
the $L^\infty$ context, we will drop $\infty$ and just write
$\CP^{ham}(Sympeo(M,\omega),id)$ for
$\CP^{ham}_\infty(Sympeo(M,\omega),id)$ from now on to simplify the
notations, unless explicitly mentioned otherwise.

\section{One parameter subgroups}
\label{sec:structure}

We consider the \emph{developing map}
$$
\operatorname{Dev}: \CP^{ham}(Symp(M,\omega),id) \to
C_m^\infty([0,1] \times M,\R):
$$
This is defined by the assignment of the normalized generating
Hamiltonian $H$ of $\lambda$, when $\lambda = \phi_H: t \mapsto
\phi_H^t$. We also consider the inclusion map
$$
\iota_{ham}: \CP^{ham}(Symp(M,\omega),id) \to \CP(Symp(M,\omega),id)
\hookrightarrow \CP(Homeo(M,\omega),id).
$$
Imitating \cite{oh:hameo1}, we call the product map
$(\iota_{ham},\operatorname{Dev})$ the \emph{unfolding map} and
denote the image thereof by \be\label{eq:QQ} \CQ_\infty:=
\operatorname{Image} (\iota_{ham},\operatorname{Dev}) \subset
\CP(Homeo(M),id) \times C^0_m([0,1] \times M,\R). \ee Then both maps
$\operatorname{Dev}$ and $\iota_{ham}$ are Lipschitz with respect to
the metric $d_{ham}$ on $\CP^{ham}(Symp(M,\omega),id)$ by definition
and so the unfolding map canonically extends to the closure
$\overline{\CQ}_\infty$ in $\CP(Homeo(M),id) \times C^0_m([0,1]
\times M,\R)$ in that we have the following continuous projections
\bea \overline{\iota}_{ham} &: \overline{\CQ}_\infty \to
\CP(Homeo(M),id) \label{eq:bariota}\\
\overline{\operatorname{Dev}} &: \overline{\CQ}_\infty \to
C^0_m([0,1]\times M,\R). \label{eq:barDev} \eea

We would like to note that by definition we also have the extension
of the evaluation map $ev_1: \CP^{ham}(Symp(M,\omega),id) \to
Symp(M,\omega) \to Homeo(M)$ to \be\label{eq:extendev1}
\overline{ev}_1: \operatorname{Image }(\overline\iota_{ham}) \to
Homeo(M). \ee

The following theorem was proved in \cite{oh:hameo1} in the
$L^{(1,\infty)}$ context. The proof for the $L^\infty$ context is
very simple which we present here.

\begin{thm}\label{thm:normal} The subset
$Hameo_\infty(M,\omega)$ is a path-connected normal subgroup of
$Sympeo(M,\omega)$ (with respect to the subspace topology).
\end{thm}
\begin{proof} We refer to \cite{oh:hameo1} for the proof of the
group property and focus on the proof of normality of the subgroup
$Hameo_\infty(M,\omega)$ in $Sympeo(M,\omega)$.
\par
Let $h \in Hameo_\infty(M,\omega)$ and $g \in Sympeo(M,\omega)$. By
definition, there exist a sequence $H_i$ of Hamiltonians such that
$\operatorname{hlim}_{i\to \infty}(\phi_{H_i},H_i) =
(\phi_{H_\infty},H_\infty)$ and $\phi_{H_\infty}^1 = h$, and a
sequence $\psi_i \in Symp(M,\omega)$ such that $\lim_{C^0}\psi_i =
g$. Obviously $\lim_{C^0} H_i\circ \psi_i = H \circ g$ and $\psi_i
\phi_{H_i}\psi_i^{-1}$ converges in $C^0$-topology and
$$
\lim_{i\to \infty}\psi_i \phi_{H_i}^1\psi_i^{-1} = ghg^{-1}.
$$
Hence by definition, we have proved $ghg^{-1} \in
Hameo_\infty(M,\omega)$ and hence the normality. Path-connectedness
is immediate since the above proof shows that any element $h \in
Hameo_{\infty}(M,\omega)$ can be connected to the identity via a
path lying in $Hameo_\infty(M,\omega)$ that is connected in
$Sympeo(M,\omega)$.
\end{proof}

In terms of this group, we also call a continuous Hamiltonian path
$\lambda$ a \emph{$C^0$-hamiltonian-continuous} map from $[0,1]$ to
$Hameo(M,\omega)$ and often denote it as $\lambda:[0,1] \to
Hameo(M,\omega)$. Similarly to the case of the interval $[0,1]$, we
can define a continuous Hamiltonian path on $[a,b]$ with $b > a$
$$
\lambda: [a,b] \to Hameo(M,\omega)
$$
to be a path such that \be\label{eq:lambda(a)} \lambda \circ
(\lambda(a))^{-1} \in \CP^{ham}_{[a,b]}(Sympeo(M,\omega),id) \ee
where we define $\CP^{ham}_{[a,b]}(Sympeo(M,\omega),id)$ similarly
as $\CP^{ham}(Sympeo(M,\omega),id)$ with $[0,1]$ replaced by
$[a,b]$.

\medskip

Next, we define
$$
\CH^0_m: = \operatorname{Image}(\overline{\operatorname{Dev}})
$$
and call any element therefrom a \emph{continuous Hamiltonian}.

We first prove the following theorem by the argument used in
\cite{oh:hameo1}.

\begin{thm}\label{Devone-one} The map $\overline{\operatorname{Dev}}:
\overline{\CQ}_\infty \to \CH^0_m([0,1]\times M,\R)$ is a bijective
map.
\end{thm}
\begin{proof}
Recalling that the map is nothing but the restriction to $
\overline{\CQ}_\infty$ of the projection
$$
\CP(Homeo(M),id) \times C^0_m([0,1]\times M,\R) \to
C^0_m([0,1]\times M,\R)
$$
and $\CH^0_m([0,1]\times M,\R)$ is defined to be its image of $
\overline{\CQ}_\infty$, $\overline{\operatorname{Dev}}$ is a
well-defined surjective map. To prove that it is also one-one, we
need to prove that if $(\lambda',H),\, (\lambda, H) \in
 \overline{\CQ}_\infty \subset \CP^{ham}(Sympeo(M,\omega),id)
\times \CH^0_m([0,1]\times M,\R)$, then $\lambda = \lambda'$.

By definition, if $(\lambda,H) \in   \overline{\CQ}_\infty$, there
exists a sequence of smooth Hamiltonians $H_i$ such that
\be\label{eq:Hi} H_i \to H \quad \mbox{in $C^0$}, \quad \overline
d(\phi_{H_i},\lambda) \to 0. \ee Applying the same argument to
$(\lambda',H)$, we obtain another sequence $H_i'$ such that
\be\label{eq:H'i} H_i^\prime\to H \quad \mbox{in $C^0$}, \quad
\overline d(\phi_{H_i^{\prime}},\lambda') \to 0. \ee Combining
(\ref{eq:Hi}) and (\ref{eq:H'i}), we in particular have
\be\label{eq:HiH'i} \|H_i - H_i'\|_\infty \to 0. \ee Now we will
prove the theorem by contradiction. Suppose $\lambda \neq \lambda'$.
Then, since we $\lambda(0) = \lambda(0) = id$, there exists $s \in
(0,1]$ such that $\lambda(s) \neq \lambda'(s)$ and so
$(\lambda(s))^{-1}\circ \lambda'(s) \neq id$. Since
$(\lambda(s))^{-1}\circ \lambda'(s)$ is a continuous map, there
exists a closed symplectic ball $B \subset M$ such that \be
(\lambda(s))^{-1}\circ \lambda'(s)(B) \cap B = \emptyset \ee Since
$B$ is compact and $(\phi_{H_i}^s)^{-1}\circ \phi_{H'_i}^s \to
(\lambda(s))^{-1}\circ \lambda'(s)$ uniformly, we have
$$
(\phi_{H_i}^s)^{-1}\circ \phi_{H'_i}^s(B) \cap B = \emptyset
$$
for all sufficiently large $i$'s.

We recall that when $H \mapsto \phi_H^1$ we have $H^s \mapsto
\phi_H^s$ where $H^s:[0,1] \times M \to \R$ is the Hamiltonian
defined by
$$
H^s(t,x) = sH(st,x).
$$
And the product $\phi_H \phi_F$ is also a Hamiltonian path which is
generated by the product Hamiltonian $H \# F$ which is defined by
$$
H \# F(t,x): = H(t,x) + F(t, \phi_H^t(x))
$$
and the inverse $\phi_H^{-1}$ is generated by the Hamiltonian
$$
\overline H (t,x): = - H(t,\phi_H(x)).
$$
Therefore we have $\overline H_i^s \# H_i'^s \mapsto
(\phi_{H_i}^s)^{-1}\circ \phi_{H'_i}^s$ where $\overline H_i^s \#
H_i'^s$ is given by \bea\label{eq:barHiH'i} \overline H_i^s \#
H_i'^s(t,x) & = & - H_i^s(t,\phi_{H_i}^{st}(x)) +
H_i'^s(t,\phi_{H_i}^{st}(x))
\nonumber \\
& = & s\left(H_i'(st,\phi_{H_i}^{st}(x)) -
H_i'(st,\phi_{H_i}^{st}(x))\right). \eea Then the energy-capacity
inequality from \cite{lal-mc1} implies \be\label{eq:HicB}
\|\overline H_i^s \# H_i'^s\|_\infty \geq \frac{1}{2}c(B) > 0 \ee
where $c(B)$ is the Gromov area of the symplectic ball $B$.  On the
other hand, it follows from \eqref{eq:barHiH'i} that for any $0 < s
\leq 1$, we have
$$
\|\overline H_i^s \# H_i'^s\|_\infty = s\|H_i - H_i'\|_\infty \leq
\|H_i - H_i'\|_\infty
$$
which converges to 0 by (\ref{eq:HiH'i}). This contradicts to
(\ref{eq:HicB}) and finishes the proof of $\lambda = \lambda'$.
\end{proof}

Combining Theorem \ref{Devone-one} with the uniqueness of
Hamiltonians, we immediately derive the following one-one
correspondence which extends the well-known correspondence between
smooth Hamiltonians and smooth Hamiltonian flows.

\begin{thm}\label{one-one} The composition map
$$
\bar{\iota}_{ham}\circ (\overline{\operatorname{Dev}})^{-1}:
\CH_m^0([0,1]\times M,\R) \to \CP^{ham}(Sympeo(M,\omega),id)
$$
provides a one-one correspondence between the two sets,
$\CP^{ham}(Sympeo(M,\omega),id)$ and $\CH^0([0,1]\times M,\R)$ under
which the following diagram commutes: \be\label{eq:commudiagram}
\begin{matrix}
C^\infty_m([0,1]\times M,\R) & \longrightarrow & \CP^{ham}(Symp(M,\omega),id)\\
\downarrow & {} & \downarrow \\
\CH_m^0([0,1]\times M,\R) &\longrightarrow
&\CP^{ham}(Sympeo(M,\omega),id)
\end{matrix}
\ee
\end{thm}
Via this correspondence, we will also denote the value $\phi_H(s)
\in Sympeo(M,\omega)$ of the continuous Hamiltonian path $\phi_H$ by
$\phi_H^s$. It is easy to check that $\phi_H^s$ itself lies in
$Hameo(M,\omega)$ since $H^s$ defined by $H^s(t,x) = sH(st,x)$ is a
continuous Hamiltonian generating $\phi_H^s$, i.e., $\phi_{H^s}(1) =
\phi_H(s) = \phi_H^s$.

The following is a natural question to ask

\begin{ques} Is $Hameo(M,\omega)$ a Lie group, or does it
contain a subgroup which is a Lie group bigger than $Ham(M)$?
\end{ques}

As a first step towards to the study of this question, we prove the
following theorem

\begin{thm}\label{autonom}
Suppose that $\lambda \in \CP^{ham}(Sympeo(M,\omega),id)$ is a
continuous Hamiltonian path and $H$ its Hamiltonian. Then the
followings are equivalent:
\begin{enumerate}
\item $\lambda$ is a one-parameter
subgroup, i.e., a path satisfying
$$
\lambda(t+s) = \lambda(t)\lambda(s).
$$
\item $H$ is time-independent, i.e., there exists
a continuous function $h:M \to \R$ such that $H = h\circ \pi (t,x)$
(i.e., $H_t = h$ everywhere in $t \in [0,1]$).
\end{enumerate}
We call any such function $h: M \to \R$ an \emph{autonomous}
continuous Hamiltonian of $(M,\omega)$ and denote by
$ham^{aut}_\infty(M,\omega)$ the set of such Hamiltonians.
\end{thm}
\begin{proof}
Suppose that $\lambda$ is a one-parameter subgroup. Let $H_i$ be a
sequence of smooth Hamiltonian functions such that $(\phi_{H_i},H_i)
\to (\lambda,H)$. We consider the identity
$$
\lambda(t+s)\lambda(s)^{-1} = \lambda(t).
$$
Then for each given $s$ we have \be\label{eq:phiitolambda}
\phi_{H_i}^{t+s}(\phi_{H_i}^s)^{-1} \to
\lambda(t+s)\lambda(s)^{-1}=\lambda(t) \ee uniformly in $t$. We
denote by $\mu_i^s = \mu_i^s(t)$ the smooth path defined by
$$
\mu_i^s(t) = \phi_{H_i}^{t+s}(\phi_{H_i}^s)^{-1}.
$$
Note that $\mu_i^s(0) = id$. A straightforward computation shows
$$
\operatorname{Dev}(\mu_i^s)(t,x) = H_i(t+s,x)
$$
i.e., the function $G_i(s)$ defined by $G_i(s)(t,x) = H_i(t+s,x)$
generates the Hamiltonian flow $t \mapsto
\phi_{H_i}^{t+s}(\phi_{H_i}^s)^{-1}$. However the latter converges
to $\lambda(t+s)\lambda(s)^{-1}=\lambda(t)$ in the metric $\overline
d$. Obviously the sequence $G_i(s)$ defines a  Cauchy sequence in
$C^0$ on $[c,d] \times M$ for any interval $[c,d]$. Then by the
uniqueness of the Hamiltonian \cite{viterbo2}, \cite{oh:locality} we
must have \be\label{eq:GinftyH} G_\infty(s) = H \ee for all $s$
where $G_\infty(s)$ is the limit
$$
G_\infty(s): = \lim_{C^0} G_i(s).
$$
On the other hand, using the convergence of $H_i \to H$ in $C^0$ we
can also write \be\label{eq:Ginfty} G_\infty(s)(t,x) = H(t+s,x) \ee
in $C^0$ for each fixed $s$. Therefore combining (\ref{eq:GinftyH})
and (\ref{eq:Ginfty}), we have proved \be\label{eq:HsH} H(\cdot + s,
\cdot) = H \ee as a $C^0$ function for all $s$. This also implies
$H_t = H_s$ for all $t, \, s$. Setting $h:M \to \R$ to be the common
function, we have proved that (1) implies (2).

Conversely suppose that $H = h$ is continuous and time-independent
i.e., $h(s+t,x) = h(t,x)$ for all $s,\, t$. We need to show
$\phi_h^{s+t} = \phi_h^t \cdot \phi_h^s$ i.e.,
\be\label{eq:autonomous} \phi_h^{s+t} \cdot (\phi_h^s)^{-1} =
\phi_h^t. \ee But the flow $t \mapsto \phi_h^{s+t} \cdot
(\phi_h^s)^{-1}$ has been shown to be generated by the Hamiltonian
$G$ defined by $G(s)(t,x) = h(t+s,x)$. By the assumption, we have
$G(s) = h$ for all $s$. Now injectivity, Theorem \ref{Devone-one},
of $\overline{\operatorname{Dev}}$ implies the flow $t \mapsto
\phi_h^{s+t} \cdot (\phi_h^s)^{-1}$ should coincide with the flow $t
\mapsto \phi_h^t$ and hence (\ref{eq:autonomous}).
\end{proof}

A similar argument gives rise to the law of \emph{conservation of
energy} under the autonomous continuous Hamiltonian flow. This
extends the well-known conservation law in the smooth autonomous
Hamiltonian flow.

\begin{thm}\label{conservation} Let $H$ be an autonomous continuous Hamiltonian and
$\phi_H$ be its flow. Then we have $H\circ \phi_H^s = H$ for all
$s$.
\end{thm}
\begin{proof} We note that the conjugate flow
$$
t \mapsto (\phi_H^s)^{-1}\phi_H^t\phi_H^s
$$
is generated by the continuous Hamiltonian $H \circ \phi_H^s$ for
each $s$: This follows by a similar argument by considering an
approximating sequence $H_i$ and recalling that the statement holds
for a smooth flow.

On the other hand, Theorem \ref{autonom} implies that $\phi_H$ is a
one-parameter subgroup and so we have
$$
(\phi_H^s)^{-1}\phi_H^t\phi_H^s = \phi_H^t
$$
which is generated by $H$. By the uniqueness result, we must have
$H\circ \phi_H^s = H$ for all $s$. This finishes the proof.
\end{proof}

\section{One-jets of continuous Hamiltonian paths}
\label{sec:one-jet}

In this section, we would like to associate a vector space to each
element $\phi \in Hameo(M,\omega)$ which would play the role of a
`tangent space' to  $Hameo(M,\omega)$ at $\phi$.

We first define the notion of `1-jets' of continuous Hamiltonian
paths at $\phi \in Hameo(M,\omega)$. For this purpose, we recall the
definition of the tangent map $\overline{\operatorname{Tan}}:
\CP^{ham}(Sympeo(M,\omega),id) \to C^0_m([0,1]\times M,\R)$ from
\cite{oh:hameo1}: this is defined by
$$
\operatorname{Tan}(\lambda)(t,x) = H(t,\phi_H^t(x))
$$
if $\lambda = \phi_H$, i.e., if $\operatorname{Dev}(\lambda) = H$.
From this, it follows that we have the identity
$$
\overline{\operatorname{Tan}}(\lambda) =
\overline{\operatorname{Dev}}(\lambda) \circ \lambda.
$$
Here we note that the right hand side
$\overline{\operatorname{Dev}}(\lambda) \circ \lambda$ is defined by
the formula
$$
(\overline{\operatorname{Dev}}(\lambda)\circ \lambda)(t,x) =
\overline{\operatorname{Dev}}(\lambda)(t,\lambda(t)(x))
$$
which we remark is well-defined as an element in $C^0([0,1]\times
M,\R)$ because $\overline{\operatorname{Dev}} (\lambda)$ is $C^0$
and $\lambda$ defines a continuous map on $[0,1]\times M$.

We extend this discussion to the set of continuous Hamiltonian paths
issued at $\phi \in Hameo(M,\omega)$.

\begin{defn} Let $\phi \in Hameo(M,\omega)$.
A continuous Hamiltonian path $\lambda$ with $\lambda(0) = \phi$ is
the $C^0$-limit of a sequence of smooth Hamiltonian path of the form
$$
t \mapsto \phi_{H_i}^t \circ \phi_i
$$
with $H_i: [0,1] \times M \to \R$ converging uniformly and $\phi_i
\to \phi$ converging in $C^0$ Hamiltonian topology.

We denote by $\CP^{ham}(Sympeo(M,\omega),\phi)$ the set of
continuous Hamiltonian path $\lambda$ with $\lambda(0) = \phi$.
\end{defn}

One can easily see that this definition is equivalent to the
existence of a sequence $(H_i,F_i)$ of pairs of Hamiltonians such
that $H_i$ and $F_i$ satisfy the following:
\begin{enumerate}
\item both $H_i, \, F_i$ converge uniformly,
\item $\phi_{F_i}$ converge uniformly and $\phi_{F_i}^1 \to \phi$ uniformly,
\item and the path $t \mapsto \phi_{H_i}^t\circ \phi_{F_i}^1$
uniformly converges to $\lambda$.
\end{enumerate}

We extend the map $\overline{\operatorname{Tan}}$ and
$\overline{\operatorname{Dev}}$ to
$\CP^{ham}(Sympeo(M,\omega),\phi)$ by defining \beastar
\overline{\operatorname{Tan}}(\lambda)(t,x)& = & H(t, \phi_H^t(\phi(x))\\
\overline{\operatorname{Dev}}(\lambda)(t,x) & = & H(t,x). \eeastar

Now we introduce the notion of 1-jets of continuous Hamiltonian
paths.

\begin{defn}\label{defn:germ} Let $h \in Hameo(M,\omega)$.
Consider two continuous Hamiltonian paths $\lambda_1, \, \lambda_2$
defined on $(-\e,\e)$ with $\lambda_1(0) = \lambda_2(0)=h$. We say
$\lambda_1 \sim \lambda_2$ at $h$ if the identity \be\label{eq:germ}
\overline{\operatorname{Tan}}(\lambda_1)(0) =
\overline{\operatorname{Tan}}(\lambda_2)(0) \quad \mbox{in $C^0(M)$}
\ee holds.
\end{defn}

The following lemma is an immediate consequence of the definition.

\begin{lem} The above relation is an equivalence relation.
\end{lem}

\begin{defn}[One-jets of continuous Hamiltonian paths]
For a given germ of continuous Hamiltonian path $\lambda$ issued at
$h$, we denote the equivalence class of $\lambda$ at $h$ by
$[\lambda]_h$ and call $[\lambda]_h$ a \emph{1-jet} of continuous
Hamiltonian paths at $h$. We denote by $\CT_h$ the set of 1-jets of
the continuous Hamiltonian paths at $h$. We denote \be\label{eq:CT}
\CT: = \bigcup_{h \in Hameo(M,\omega)} \CT_h. \ee
\end{defn}

We now equip $\CT_h$ with a vector space structure. For the
addition, we start with the following lemma whose proof is immediate
from the definitions.

To describe $\CT_h$ more concretely, let us first consider the case
$h=id$.

\begin{lem}\label{ev0Dev}
\begin{enumerate}
\item
Let $\lambda_1, \, \lambda_2 \in \CP^{ham} (Sympeo(M,\omega),id)$
satisfy $\lambda_1(0) = \lambda_2(0) =id$. Denote by $ev_0$ the
evaluation at $t = 0$, i.e., $ev_0(H) = H_0$. Let
$h=ev_0(\overline{\operatorname{Dev}}(\lambda_1))$ and $g =
ev_0(\overline{\operatorname{Dev}}(\lambda_2))$. Then $h+g$ lies in
the image of $ev_0\circ \overline{\operatorname{Dev}}$.
\item Let $\lambda \in \CP^{ham} (Sympeo(M,\omega),id)$
and $k = \overline{\operatorname{Dev}}(\lambda)(0)$ and $a \in \R$.
Then $a \cdot k$ lies in the image of $ev_0\circ
\overline{\operatorname{Dev}}$.
\end{enumerate}
\end{lem}
\begin{proof} We consider the product path
$\lambda = \lambda_1\cdot \lambda_2$. Obviously $\lambda(0) =
\lambda_1(0)\cdot \lambda_2(0) = id$. On the other hand, the
identity \be\label{eq:bar-product}
\overline{\operatorname{Dev}}(\lambda) =
\overline{\operatorname{Dev}}(\lambda_1) +
\overline{\operatorname{Dev}}(\lambda_2)\circ \lambda_1^{-1} \ee as
a $C^0$-function was proved in the proof of Theorem 3.23, especially
(3.28) of \cite{oh:hameo1}. Recall this formula for smooth
Hamiltonian paths is well-known. We evaluate the above identity at
$t = 0$, which first implies that $\lambda = \lambda_1\cdot
\lambda_2$ has its value given by
$$
\overline{\operatorname{Dev}}(\lambda)(0) =
\overline{\operatorname{Dev}}(\lambda_1)(0)
+\overline{\operatorname{Dev}}(\lambda_2)(0) = h+g.
$$

For the second statement with $a \neq 0$, we consider the rescaled
path $\lambda^a$ defined by
$$
\lambda^a(t) = \lambda(at).
$$
Then it follows that
$\overline{\operatorname{Dev}}(\lambda^a)(t,\cdot) = a
\overline{\operatorname{Dev}}(\lambda)(at,\cdot)$. This identity
shows that $\lambda^a$ has the value given by
$\overline{\operatorname{Dev}}(\lambda^a)(0) =
a\overline{\operatorname{Dev}}(\lambda)(0) = a h$ which proves the
statement for $a \neq 0$. When $a = 0$, we just consider the
constant path $\lambda \equiv id$. This finishes the proof.
\end{proof}

For two given continuous Hamiltonian paths $\lambda,  \lambda'$
issued at $id$, Theorem \ref{Devone-one} and the above lemma enable
us to define the sum \be\label{eq:addition} [\lambda_1]_{id} +
[\lambda_2]_{id} = [\lambda_1\cdot\lambda_2]_{id}. \ee For the
scalar multiplication, we define \be\label{eq:scalarmult} a \cdot
[\lambda]_{id} = [\lambda^a]_{id}. \ee

For the general continuous Hamiltonian paths $\lambda_1, \,
\lambda_2$ at $\phi \in Hameo(M,\omega)$, we note that
$\lambda_1\phi^{-1}\lambda_2$ is a continuous Hamiltonian path at
$\phi$ and so we define
$$
[\lambda_1]_\phi + [\lambda_2]_\phi =
[\lambda_1\phi^{-1}\lambda_2]_\phi.
$$
Similarly we define
$$
a \cdot [\lambda]_\phi = [\lambda^a]_\phi
$$
noting that $\lambda^a$ is a continuous Hamiltonian path with
$\lambda^a(0) = \phi$ if $\lambda(0) = \phi$.

The following is straightforward to check from the definition.

\begin{prop} The set $\CT_{\phi}$ forms a vector space. And the maps
$$
\overline{\operatorname{Tan}}_\phi, \,
\overline{\operatorname{Dev}}_\phi: \CT_{\phi} \to  C^0(M),
$$
defined by
$$
\overline{\operatorname{Tan}}_\phi([\lambda]_\phi) =
\overline{\operatorname{Tan}}(\lambda)(0)
$$
and
$$
\overline{\operatorname{Dev}}_\phi([\lambda]_\phi) =
\overline{\operatorname{Dev}}(\lambda)(0)
$$
respectively, define injective homomorphisms. Furthermore the image
of $\overline{\operatorname{Dev}}_\phi$ is independent of $\phi \in
Hameo(M,\omega)$.
\end{prop}
\begin{proof} We will just show the identity
$$
a\cdot ([\lambda_1]_\phi + [\lambda_2]_\phi) = a\cdot
[\lambda_1]_\phi + a\cdot [\lambda_2]_\phi
$$
and leave the rest to the readers.

We compute
$$
a\cdot ([\lambda_1]_\phi + [\lambda_2]_\phi) = a \cdot [\lambda_1
\phi^{-1} \lambda_2]_\phi = [\lambda_1^a \phi^{-1} \lambda_2^a]_\phi
$$
and
$$
a\cdot [\lambda_1]_\phi + a\cdot [\lambda_2]_\phi =
[\lambda^a_1]_\phi + [\lambda^a_2]_\phi = [\lambda^a_1 \phi^{-1}
\lambda^a_2]_\phi.
$$
For the last identity, we use the fact $\lambda^a_1(0) =
\lambda^a_2(0) = \phi$ for all $a$. This finishes the proof.
\end{proof}

Now we denote by $ham_\infty(M,\omega)$ the subset
\be\label{eq:ham1} ham_\infty(M,\omega): = \operatorname{Im
}(\overline{\operatorname{Dev}}_{id}) \subset C^0(M). \ee Then Lemma
\ref{ev0Dev} and the above discussion imply that
$ham_\infty(M,\omega)$ is a subspace of the vector space $C^0(M)$.

Therefore the union $\CT \to Hameo(M,\omega)$ forms a `vector
bundle' with a canonical trivialization
$$
\overline{\operatorname{Dev}}: \CT \to Hameo(M,\omega) \times
ham_\infty(M,\omega)
$$
such that the following diagram commutes: \be\label{eq:Tham}
\begin{matrix}
T(Ham(M,\omega)) & \longrightarrow & Ham(M,\omega) \times C^\infty_m(M) \\
\downarrow & {} & \downarrow \\
\CT & \longrightarrow & Hameo(M,\omega) \times ham_\infty(M,\omega)
\end{matrix}
\ee where the horizontal maps are induced by the developing map
$\overline{\operatorname{Dev}}$ defined above.
\begin{defn}
We call $\CT \to Hameo(M,\omega)$ the \emph{hamiltonian tangent
bundle} of $Hameo(M,\omega)$.
\end{defn}

\begin{defn} Let $\lambda:[a,b] \to Hameo(M,\omega)$ be
a continuous Hamiltonian path. We denote
$$
\lambda'(s): = [\lambda]_{\lambda(s)}
$$
and call the \emph{hamiltonian tangent vector} of the path.
\end{defn}
By definition of the equivalence class $[\lambda]_{\phi}$ in
Definition \ref{defn:germ} we can identify $\lambda'(s)$ with $
\overline{\operatorname{Tan}}(\lambda)(s). $ Under this
identification, we also have
$$
\lambda'(s) \circ \lambda(s)^{-1} =
\overline{\operatorname{Dev}}(\lambda)(s).
$$
The above discussion somehow indicates that all continuous
Hamiltonian path `differentiable' and so carries a `tangent vector
field' which is continuous.

The following questions seem to be important questions to ask.

\begin{ques}
\begin{enumerate}
\item It is easy to see from definition that
\be\label{eq:inclusion} ham_\infty^{aut}(M,\omega)  \subset
ham_\infty(M,\omega) \subset C^0(M). \ee Are any of these inclusions
strict?
\item Is $ham_\infty^{aut}(M,\omega)$ a subspace of $C^0(M)$ too?
\item Do we have a `smooth' structure on $Hameo(M,\omega)$ so that
$\CT$ becomes its tangent bundle? \item Is $Hameo(M,\omega)$ a Lie
group? \item Can we define a Poisson bracket on
$ham_\infty(M,\omega)$? Or what is the maximal subspace of
$ham_\infty(M,\omega)$ on which the bracket operation is defined in
a way that it extends the standard Poisson bracket on
$C^\infty_m(M)$?
\end{enumerate}
\end{ques}

We believe that the following conjecture is true.

\begin{conj} The subset $ham_\infty(M,\omega)$ is a proper subset of $C^0(M)$.
\end{conj}

We like to compare this conjecture to the following group analog
proposed in \cite{oh:hameo1}

\begin{conj}[\cite{oh:hameo1}] $Hameo(M,\omega)$ is a proper
subgroup of $Sympeo_0(M,\omega)$ in general.
\end{conj}

It was shown in \cite{oh:hameo1} that this conjecture is true
whenever the mass flow homomorphism is non-trivial or there exists a
symplectic diffeomorphism that has no fixed point, e.g., $T^{2n}$.

\section{Extended Hofer length and the intrinsic norm}
\label{sec:hoferlength}

In this section, using the uniqueness result of continuous
Hamiltonians, we will extend the definition of the Hofer length
function to the continuous Hamiltonian paths, and define its
associated intrinsic distance function on $Hameo(M,\omega)$.

First, the uniqueness theorem of Hamiltonian of a continuous
Hamiltonian path enables us to define the following extension of the
Hofer length to continuous Hamiltonian paths.

\begin{defn}\label{length}
Let $\lambda \in \CP^{ham}(Sympeo(M,\omega), id)$. We define the
length of a continuous Hamiltonian path
$$
\lambda: \CP^{ham}(Sympeo(M,\omega), id)  \to \R_+
$$
by \be\label{eq:leng} \operatorname{leng}(\lambda): =
\|\overline{\operatorname{Dev}}(\lambda)\| = \int_0^1 \osc H_t\, dt
\ee where $\overline{\operatorname{Dev}}(\lambda) = H$. We call this
the {\it Hofer length} of the continuous Hamiltonian path $\lambda$.
\end{defn}

The uniqueness theorem implies that $\operatorname{leng}(\lambda)$
is also the same as the limit \be\label{eq:lengconv} \lim_{i \to
\infty} \operatorname{leng}(\phi_{H_i}) = \lim_{i \to \infty}
\|H_i\| \ee for any sequence $(\phi_i,H_i) \to (\lambda,H)$. In
particular, the definition extends that of the smooth case.

\begin{prop}\label{hoferlength} The function
$\operatorname{leng}: \CP^{ham}(Sympeo(M,\omega),id) \to \R$
satisfies the triangle inequality \be\label{eq:leng-tri}
\operatorname{leng}(\lambda\mu) \leq \operatorname{leng}(\lambda) +
\operatorname{leng}(\mu) \ee and is continuous.
\end{prop}
\begin{proof} The triangle inequality
(\ref{eq:leng-tri}) is an immediate consequence of
(\ref{eq:bar-product}), i.e.,
$$
\overline{\operatorname{Dev}}(\lambda\mu) =
\overline{\operatorname{Dev}}(\lambda) +
\overline{\operatorname{Dev}}(\mu) \circ \lambda^{-1}.
$$
The triangle inequality then gives rise to the inequality
$$
|\operatorname{leng}(\lambda) - \operatorname{leng}(\mu)| \leq
\operatorname{leng}(\lambda^{-1}\mu)
$$
from which continuity of $\operatorname{leng}$ follows.
\end{proof}

Next we consider the action of $Sympeo(M)$ on
$$
\CP(Homeo(M),id) \times C^0([0,1] \times M)
$$
given by \be\label{eq:sympeoact} \left(\psi, (\lambda,H)\right)
\mapsto (\psi^{-1}\lambda \psi, H\circ \psi) \ee and we prove the
invariance property of the length under this action. In
\cite{oh:hameo1}, the action (\ref{eq:sympeoact}) was proven to map
$\CP^{ham}(Sympeo(M,\omega),id)$ to itself and so induce an action
thereon. We denote this action on $\CP^{ham}(Sympeo(M,\omega),id)$
by
$$
\psi\cdot \lambda = \psi^{-1}\lambda \psi.
$$

Now we prove the following theorem
\begin{thm}\label{lenginvariant} The length function defined in (\ref{eq:leng}) is
invariant under the action (\ref{eq:sympeoact}) on
$\CP^{ham}(Sympeo(M,\omega),id)$.
\end{thm}
\begin{proof}
Let $\psi \in Sympeo(M,\omega)$ and $\lambda \in
\CP^{ham}(Sympeo(M,\omega),id)$. Then
$$
\psi \cdot \lambda = \lim_{i \to \infty}\psi_i^{-1}\phi_{H_i}\psi_i
$$
with $H_i\circ \psi_i$ converging for any sequence $\psi_i \in
Symp(M,\omega)$ with $\psi_i \to \psi$ uniformly and
$(\phi_{H_i},H_i) \to (\lambda,H_\infty)$ converging in
$C^0$-Hamiltonian topology. Therefore from the definition
(\ref{eq:leng}) and convergence (\ref{eq:lengconv}), we derive
$$
\operatorname{leng}(\psi\cdot \lambda) = \lim_{i \to \infty}
\operatorname{leng}(\psi_i^{-1}\phi_{H_i}\psi_i) = \lim_{i\to
\infty}\operatorname{leng}(\phi_{H_i}) =
\operatorname{leng}(\lambda)
$$
which finishes the proof.
\end{proof}

Next we recall the definition of Hofer displacement energy $e(A)$:
for every compact subset $A \subset M$,
$$
e(A):= \inf_{H}\{ \|H \| \mid H \in C^\infty([0,1] \times M,\R) , \,
A \cap \phi_H^1(A) = \emptyset \}
$$
We can generalize this generalized context of continuous
Hamiltonians: For every compact subset $A \subset M$, we define
$$
\overline e(A):= \inf_{\lambda}\{ \leng(\lambda) \mid \lambda \in
\CP^{ham}(Sympeo(M,\omega),id), \, A \cap \lambda(1)(A) = \emptyset
\}.
$$
Obviously we have $e(A) \geq \overline e(A)$. In addition, we prove
\begin{thm}\label{e=bare}
For every compact $A \subset M$, we have
$$
e(A) = \overline e(A).
$$
\end{thm}
\begin{proof}
For the opposite inequality, let $\delta > 0$. By definition of
$\overline e(A)$, we have $\lambda \in
\CP^{ham}(Sympeo(M,\omega),id)$ such that $\lambda(1)(A) \cap A =
\emptyset$ and
$$
\leng(\lambda) \leq \overline e(A) + \frac{\delta}{2}.
$$
Let $H_i$ be a sequence of smooth Hamiltonians with
$\operatorname{hlim}(\phi_{H_i}) = \lambda$. Then for all
sufficiently large $i$ we have
$$
\|H_i\| \leq  \operatorname{leng}(\lambda) + \frac{\delta}{2} \leq
\overline e(A) + \delta.
$$
On the other hand, we have
$$
e(A) \leq \|H_i\|
$$
since we have $\phi_{H_i}^1(A) \cap A = \emptyset$ for all
sufficiently large $i$ by the $C^0$-convergence of $\phi_{H_i}^1 \to
h$, $h(A) \cap A = \emptyset$ and compactness of $A$. Altogether we
have derived
$$
e(A) \leq \overline e(A) + \delta.
$$
Since $\delta > 0$ is arbitrary, this implies $e(A) \leq \overline
e(A)$. This finishes the proof.
\end{proof}

Based on this theorem, we will just denote by $e(A)$ the Hofer
displacement energy even in the continuous Hamiltonian category. The
following is an analog to the well-known fact that $e$ is invariant
under the action of symplectic diffeomorphisms whose proof we omit
referring to the proof of the more non-trivial case of spectral
displacement energy in the next section.

\begin{cor}\label{egA=eA} Let $g \in Sympeo(M,\omega)$ and $A$ be a compact subset of
$M$. Then we have
$$
e(g(A)) =  e(A).
$$
\end{cor}

Now we can use the Hofer-length function generalized to the set of
continuous Hamiltonian paths, and define an intrinsic norm of
Hofer-type on $Hameo(M,\omega)$ which in turn induces a bi-invariant
distance on $Hameo(M,\omega)$.

\begin{defn}[Intrinsic norm]\label{norm} For any $h \in Hameo(M,\omega)$,
we define
$$
\|h\| = \inf_{\lambda} \{\, \mbox{\rm leng}(\lambda) \mid \lambda
\in \CP^{ham}(Sympeo(M,\omega), id), \, \overline{ev}_1(\lambda) = h
\}.
$$
Then we define an invariant distance function
$$
d: Hameo(M,\omega) \times Hameo(M,\omega) \to \R_+
$$
by $d(h,k) = \|h^{-1}k\|$.
\end{defn}

The following theorem is the continuous Hamiltonian analog to the
well-known theorem on the Hofer norm on $Ham(M,\omega)$
\cite{hofer2}.

\begin{thm}\label{extendedHofer} Let $g, \, h \in Hameo(M,\omega)$. Then
the extended Hofer norm function
$$
\| \cdot \|: Hameo(M,\omega) \to \R_+
$$
is continuous in the Hamiltonian topology, and satisfies the
following properties:
\begin{enumerate}
\item (Symmetry) \hskip0.72in $\|g\| = \|g^{-1} \|$

\item (Triangle inequality) \hskip 0.2in $\|gh\| \leq \|g\| + \|h\|$

\item (Symplectic invariance) $\|\psi^{-1}g \psi\| = \|g\|$ for
any $\psi \in Sympeo(M,\omega)$,
\item (Nondegeneracy) \hskip0.5in $g = id$ if and only if $\|g\| = 0$.
\end{enumerate}
\end{thm}
\begin{proof} The continuity is immediate from that of
the length function $\operatorname{leng}$ in Definition \ref{length}
and from the definition of $Hameo(M,\omega)$ in (\ref{eq:hameo}).
The symmetry is straightforward to check.

For the symplectic invariance, we need to prove the identity $\|g\|
= \|\psi^{-1} g \psi\|$. According to Definition 4.2 (2), we have
\be\label{eq:normpsig} \|\psi^{-1} g \psi\| \leq
\inf_{\lambda}\left\{ \operatorname{leng}(\psi^{-1}\lambda\psi) \mid
\lambda \in \CP^{ham}(Sympeo(M,\omega),id), ev_1(\lambda) = g
\right\} \ee since it follows that if $\lambda(1) = h$, then
$\psi^{-1}\lambda(1)\psi = \psi^{-1}h \psi$. By the invariance
Theorem \ref{lenginvariant}, we have
$\operatorname{length}(\psi^{-1}\lambda\psi) =
\operatorname{length}(\lambda)$. Substituting this into
(\ref{eq:normpsig}), we have proven $\|\psi^{-1} g \psi\| \leq
\|g\|$. Applying the same argument with $\psi$ replaced by
$\psi^{-1}$ and $g$ replaced by $\psi^{-1}g\psi$, we have also
obtained $\|\psi^{-1} g \psi\| \geq \|g\|$, which finishes the
proof.

Next we prove the triangle inequality and nondegeneracy in detail.
Let $\delta > 0$ be given. By Definition \ref{norm}, there exist
$\lambda, \, \mu \in \CP^{ham}(Sympeo(M,\omega),id)$ such that
$$
\overline{ev}_1(\lambda) = g, \quad \overline{ev}_1(\mu) = h
$$
and
$$
\operatorname{leng}(\lambda) \leq \|g\| + \frac{\delta}{2}, \quad
\operatorname{leng}(\mu) \leq \|h\| + \frac{\delta}{2}.
$$
By definition, there exist sequences $H_i, \, F_i$ of smooth
Hamiltonians such that
$$
\operatorname{hlim}_{i\to \infty}(\phi_{H_i},H_i) = \lambda, \quad
\operatorname{hlim}_{i\to \infty}(\phi_{F_i},F_i) = \mu.
$$
On the other hand, since $\overline{ev}_1(\lambda\mu) = gh$, we have
$$
\|gh\| \leq \operatorname{leng}(\lambda\mu) \leq
\operatorname{leng}(\lambda) + \operatorname{leng}(\mu) \leq  \|g\|
+ \|h\| + \delta.
$$
Since $\delta$ is arbitrary, we have proven the triangle inequality.

Finally we prove nondegeneracy. Suppose that $id \neq g \in
Hameo(M,\omega)$. Since $g \neq id$ is a homeomorphism, there exists
a small symplectic ball $B(u)$ such that $g(B(u)) \cap B(u) =
\emptyset$. Let $\lambda \in \CP^{ham}(Sympeo(M,\omega),id)$ be any
given element with $\overline ev_1(\lambda) = g$. Choose a sequence
$(\phi_i,H_i)$ such that $\operatorname{hlim}(\phi_i,H_i) =
\lambda$. Then we have \be\label{eq:limHi} \lim_{i \to
\infty}\|H_i\| = \operatorname{leng}(\lambda) \ee by
(\ref{eq:lengconv}). Since $B(u)$ is compact and $g(B(u)) \cap B(u)
= \emptyset$, we also have \be \label{eq:phiBu} \phi_i(B(u)) \cap
B(u) = \emptyset \ee for all sufficiently large $i$ because $\phi_i
\to g$ in $C^0$ topology. By definition of the Hofer displacement
energy, (\ref{eq:phiBu}) implies \be\label{eq:e(B(u))} \|H_i\| \geq
\|\phi_i\| \geq e(B(u))>0 \ee for all sufficiently large $i$. The
latter positivity follows from the energy-capacity inequality proven
in [LM].

Then (\ref{eq:limHi}) and (\ref{eq:e(B(u))}) imply
$\operatorname{leng}(\lambda) \geq e(B(u)) > 0$. Since this is true
for any $\lambda \in\CP^{ham}(Sympeo(M,\omega),id)$ with
$\overline{ev}_1(\lambda) = g$, this gives rise to $\|g\| \geq
e(B(u)) > 0$ and finishes the proof of nondegeneracy.
\end{proof}

However Stefan M\"uller \cite{mueller} pointed out that the answer
to the following question is open

\begin{ques}[M\"uller]\label{mueller} Denote by $\|\cdot \|_{Ham}$ and
$\|\cdot\|_{Hameo}$ the Hofer norm on $Ham(M,\omega)$ and the
extended Hofer norm on $Hameo(M,\omega)$. Let $\phi \in
Ham(M,\omega) \subset Hameo(M,\omega)$. Does the following identity
$$
\|\phi\|_{Ham} = \|\phi\|_{Hameo}
$$
hold in general?
\end{ques}

This question can be shed some light on by relating it to the
general construction of \emph{path metric spaces} $(X,d_\ell)$
starting from a general metric space $(X,d)$ in the point of view of
Chapter 1 \cite{gromov:metric}, although the question is not exactly
in the context of this general construction therein because of its
interplay with the $C^0$-metric in addition.

\section{Spectral invariants of continuous Hamiltonian paths}
\label{sec:spectral}

In this section, we extend the definition and basic properties of
the spectral invariants of \emph{Hamiltonian paths} formulated in
\cite{oh:alan} to continuous Hamiltonian category.

For this extension, it is crucial to have the definition in the
level of Hamiltonian paths, i.e., on $\CP^{ham}(Symp(M,\omega),id)$
as formulated in \cite{oh:alan}, not just on the covering space of
$Ham(M,\omega)$. We refer to \cite{viterbo1}, \cite{oh:jdg,oh:cag},
\cite{mschwarz} for the earlier definition of similar invariants in
the context of \emph{exact} cases. Furthermore the uniqueness
theorem of continuous Hamiltonians will be crucial for the extension
to the $C^0$ category.

We first recall the definition and basic properties of the spectral
invariants $\rho(H;a)$ for a \emph{time-periodic} Hamiltonian $H$
from \cite{oh:alan}, but with some twists to incorporate the
Hamiltonian topology in its presentation.

For a given time-periodic Hamiltonian and a choice of time-periodic
almost complex structure $J$, we consider the perturbed
Cauchy-Riemann equation
$$
\frac{\partial u}{\partial \tau} + J \left(\frac{\partial
u}{\partial t} - X_{H}(t,u) \right) = 0
$$
and its associated Floer complex $\partial_{(H,J)}: CF(H) \to
CF(H)$. $CF(H)$ has a natural (increasing) filtration $CF^\lambda(H)
\hookrightarrow CF(H)$ for $\lambda \in \R$ induced by the action
functional $\CA_H$.

We call a Floer chain $\alpha$ a \emph{Floer cycle} if
$\partial_{(H,J)} \alpha = 0$ and denote its homology class by
$[\alpha]$. We define and denote the \emph{level} of the chain
$\alpha$ by
$$
\lambda_H(\alpha) = \max \{ \lambda \mid \alpha \in CF^{\lambda}(H)
\}.
$$

\begin{defn}[Definition \& Theorem 7.7 \cite{oh:alan}]
Let $H$ be a time-periodic Hamiltonian. Let $a \neq 0$ be a given
quantum cohomology class in $QH^*(M)$, and denote by $a^\flat \in
FH_*$ the Floer homology class dual to $a$ in the sense of
\cite{oh:alan}. For any given Hamiltonian path $\lambda = \phi_H \in
\CP^{ham}(Symp(M,\omega), id)$ such that $H$ is non-degenerate in
the Floer theoretic sense, we define
$$
\rho(\lambda;a): = \rho(H;a) = \inf_{\alpha \in \ker \partial_H} \{
\lambda_H(\alpha) \mid [\alpha] = a^\flat \}
$$
where $a^\flat$ is the dual to the quantum cohomology class $a$ in
the sense of \cite{oh:alan}. Then this number is finite for any
quantum cohomology class $a \neq 0$. We call any of these {\it
spectral invariants} of the Hamiltonian path $\lambda$.
\end{defn}
We refer readers to \cite{oh:alan} for the complete discussion on
general properties of $\rho(H;a)$.

Now let $H: [0,1] \times M \to \R$ be any smooth Hamiltonian,
\emph{not necessarily periodic} and let $\lambda = \phi_H$ be its
Hamiltonian path. We now explain how we associate the spectral
invariant $\rho(\lambda;a)$ to such a path $\lambda$.

Out of the given Hamiltonian $H$, we consider the time-periodic
Hamiltonian of the type $H^\zeta$ where $\zeta$ is a
reparameterization of $[0,1]$ of the type \be\label{eq:zeta} \zeta
(t) = \begin{cases}
             0  & \mbox{for }\, 0 \leq t \leq \frac{\e_0}{2}\\
             1  & \mbox{for }\, 1 - \frac{\e_0}{2} \leq t \leq 1
             \end{cases}
\ee and
$$
\zeta'(t) \geq 0 \quad \mbox{for all} \quad t \in [0,1],
$$
and the reparameterized Hamiltonian by $H^\zeta$ is given by
$$
H^{\zeta} (t,x) = \zeta'(t)H(\zeta(t),x)
$$
which generates the Hamiltonian isotopy $ t \mapsto
\phi_H^{\zeta(t)}$ in general. The following norm,
\be\label{eq:ham-norm} \|\zeta -id\|_{ham} = \|\zeta - id\|_{C^0} +
\int_0^1 |\zeta'(t) - 1|\, dt. \ee which measures a distance of
$\zeta$ from the identity parametrization, turns out to be useful as
illustrated in \cite{oh:hameo1}.

To assign a well-defined number $\rho(\lambda;a)$ depending only on
$H$ itself not on its reparameterization $H^\zeta$, we note that any
two such reparameterized Hamiltonian paths are homotopic to each
other. The homotopy invariance axiom of the spectral invariants from
\cite{oh:alan,oh:minimax} imply that the following definition is
well-defined in that it does not depend on the choice of $\zeta$.

\begin{defn} Let $\lambda$ be any smooth Hamiltonian path and
$H$ be its generating Hamiltonian. We pick a  $\zeta$ so that
$\|\zeta - id\|_{ham}$ so small that all the properties in the
$C^0$-approximation Lemma \cite{oh:hameo1} hold. Then we define
\be\label{eq:defnrho} \rho(\lambda;a):= \rho(H^\zeta;a). \ee
\end{defn}

In \cite{oh:ajm1}, \cite{oh:alan}, we proved the general inequality
\be\label{eq:maxrhomin} \int_0^1 -\max_x(K-H)\,dt \leq \rho(K;a) -
\rho(H;a) \leq \int_0^1 - \min_x(K-H)\, dt \ee for two nondegenerate
Hamiltonian functions $H, \, K$. This enabled us to extend the
definition of $\rho(\cdot; a)$ to arbitrary smooth Hamiltonian $H$
by setting
$$
\rho(H;a) = \rho(H^\zeta;a).
$$

\begin{thm}\label{spectral} For a smooth Hamiltonian path $\phi_H$, we define
$$
\rho(\phi_H;a) =  \rho(H;a).
$$
Then the map $\rho_a: \phi_H \mapsto \rho(\phi_H;a)$ extends to a
continuous function
$$
\overline\rho_a = \overline\rho(\cdot; a):
\CP^{ham}(Sympeo(M,\omega),id) \to \R
$$
(in the Hamiltonian topology) and satisfies the triangle inequality
\be\label{eq:tri-rho} \overline\rho(\lambda\mu;a\cdot b) \leq
\overline\rho(\lambda;a) + \overline\rho(\mu;b). \ee
\end{thm}
\begin{proof}
The first statement is an immediate consequence of Hamiltonian
continuity of
$$
\rho_a: \CP^{ham}(Symp(M,\omega),id) \to \R
$$
and the uniqueness theorem  of continuous Hamiltonians from
\cite{viterbo2}, \cite{oh:locality}: for any continuous Hamiltonian
path $\lambda$, we define \be\label{eq:rholambdaa}
\overline\rho(\lambda;a) = \lim_{i \to \infty}\rho(H_i;a) \ee for
any Cauchy sequence $(\phi_{H_i},H_i) \to \lambda$. The uniqueness
theorem implies that this definition is well-defined. And then
(\ref{eq:maxrhomin}) proves continuity of the extension on
$\CP^{ham}(Sympeo(M,\omega),id)$.

For the proof of triangle inequality, choose any smooth sequences
$\lambda_i$ and $\mu_i$ converging to $\lambda$ and $\mu$
respectively in the Hamiltonian topology. For smooth Hamiltonian
paths, the inequality \be \label{eq:triangle}
\rho(\lambda_i\mu_i;a\cdot b) \leq \rho(\lambda_i;a) +
\rho(\mu_i;b). \ee was proven in \cite{oh:alan} (See \cite{mschwarz}
also for the exact case). Using the continuity of $\rho$ and taking
the limit of this inequality, we have proved (\ref{eq:tri-rho}).
\end{proof}

Now we focus on the invariant $\rho(\lambda;1)$ for $1 \in QH^*(M)$.
Recall the function
$$
\gamma(H) = \rho(H;1) + \rho(\overline H;1)
$$
was introduced for a smooth Hamiltonian path $\lambda = \phi_H$ in
\cite{oh:alan,oh:dmj}. We will change its notation here to
$\operatorname{norm}_\gamma(\lambda)$ not to confuse it with the
same notation used for the spectral norm function $\gamma:
Ham(M,\omega) \to \R$ below. The function
$\operatorname{norm}_\gamma$ was proven to be non-negative and to
depend only on the path-homotopy class of $\lambda=\phi_H$ in
$Ham(M,\omega)$.

\begin{defn}[Spectral pseudo-norm] Let $\lambda \in \CP^{ham}(Symp(M,\omega),id)$ and $H$
be a Hamiltonian such that $\lambda = \phi_H$. Then we define the
function
$$
\operatorname{norm}_\gamma: \CP^{ham}(Symp(M,\omega),id) \to \R_+
$$
by setting $\operatorname{norm}_\gamma(\lambda) = \gamma(H)$. We
call $\operatorname{norm}_\gamma(\lambda)$ the \emph{spectral
pseudo-norm} of $\lambda$.
\end{defn}

Again the uniqueness of continuous Hamiltonians enables us to extend
the definition to the continuous Hamiltonian paths.

\begin{prop}\label{extlambda} The spectral pseudo-norm function $\operatorname{norm}_\gamma$ extends to a
continuous function
$$
\operatorname{norm}_{\overline\gamma}:
\CP^{ham}(Sympeo(M,\omega),id) \to \R
$$
with the definition \be\label{eq:extlambda}
\operatorname{norm}_{\overline\gamma}(\lambda) = \lim_{i\to
\infty}\operatorname{norm}_\gamma(\phi_{H_i}) \ee for a (and so any)
sequence $(\phi_i,H_i) \to \lambda$ in $C^0$-Hamiltonian topology.
\end{prop}

\begin{proof}
The proof is similar to that of the Hofer length
$\operatorname{leng}$ in that it is based on the uniqueness of
Hamiltonians and the triangle inequality
$$
\operatorname{norm}_\gamma(\lambda\mu) \leq
\operatorname{norm}_\gamma(\lambda) +
\operatorname{norm}_\gamma(\mu)
$$
and so omitted.
\end{proof}

Recall that for a smooth Hamiltonian $H$ each
$\rho(\phi_H;a)=\rho(H;a)$ is associated to a periodic orbit of
Hamilton's equation $\dot x = X_H(t,x)$ and corresponds to the
action of the periodic orbit, at least for the rational symplectic
manifold. (See \cite{oh:alan,oh:minimax}.) In this regard, the
following question seems to be of fundamental importance.

\begin{ques} What is the meaning of the extended spectral pseudo-norm $\operatorname{norm}_{\overline\gamma}(\lambda)$  in regard
to the dynamics of continuous Hamiltonian flows?
\end{ques}

In \cite{oh:dmj}, the author has introduced the notion of {\it
spectral displacement energy}. The following is the analog of the
definition from \cite{oh:dmj} of the spectral displacement energy in
the continuous Hamiltonian category.

\begin{defn}[Spectral displacement energy]
Let $A \subset M$ be a compact subset. We define the spectral
displacement energy, denoted by $e_{\overline \gamma}(A)$, of $A$ by
$$
e_{\overline \gamma}(A) =
\inf_{\lambda}\{\operatorname{norm}_{\overline\gamma}(\lambda) \mid
A \cap \lambda(1)(A) = \emptyset, \, \lambda \in
\CP^{ham}(Sympeo(M,\omega),id) \}.
$$
\end{defn}

By unraveling the definitions of Hamiltonian homeomorphisms and of
the spectral displacement energy, we also have the following theorem
whose proof will be the same as the Hofer displacement energy case
and so omitted.

\begin{thm} We have $e_{\overline \gamma}(A) = e_\gamma(A)$
for any compact subset $A \subset M$.
\end{thm}

Again based on this theorem, we just denote the spectral
displacement energy of $A$ even in the continuous Hamiltonian
category by $e_\gamma(A)$. Then we have the following theorem

\begin{thm}
For every $\psi \in Sympeo(M,\omega)$ we have
$$
e_\gamma(A) = e_\gamma(\psi(A)).
$$
\end{thm}
\begin{proof} We note that
$h(A) \cap A = \emptyset$ if and only if $\psi h\psi^{-1}(\psi(A))
\cap \psi(A) = \emptyset$. Furthermore $h \in Hameo(M,\omega)$ if
and only if $\psi h\psi^{-1}\in Hameo(M,\omega)$. This combined with
the conjugation invariance of the Hofer length finishes the proof.
\end{proof}

Next we recall that in \cite{oh:alan} we introduced the non-negative
function
$$
\gamma(\phi):= \inf_{ev_1(\lambda) =
\phi}\operatorname{norm}_{\gamma}(\lambda) = \inf_{H \mapsto
\phi}\{\rho(H;1) + \rho(\overline H;1)\}.
$$
and proved that it satisfies the properties of a bi-invariant norm
on $Ham(M,\omega)$ which we called the {\it spectral norm}. The
following definition extends this definition to $Hameo(M,\omega)$.

\begin{defn}[Spectral norm]
Let $h \in Hameo(M,\omega)$ and consider continuous Hamiltonian
paths $\lambda \in \CP^{ham}(Sympeo(M,\omega),id)$ with
$\overline{ev}_1(0) = h$. We denote by $\lambda \mapsto h$ if
$\overline{ev}_1(\lambda) = h$. We define $\overline\gamma$ by
\be\label{eq:extgamma} \overline{\gamma}(h) = \inf_{\lambda}
\{\operatorname{norm}_{\overline\gamma}(\lambda) \mid \lambda \in
\CP^{ham}(Sympeo(M,\omega),id), \, \overline{ev}_1(\lambda) = h\}.
\ee
\end{defn}

The following establishes the analogs to all the properties of
invariant norm in this continuous Hamiltonian context.

\begin{thm}\label{extendedgamma} The generalized spectral
function $\overline \gamma: Hameo(M,\omega) \to \R_+$ satisfies all
the properties of an invariant norm stated in Theorem
\ref{extendedHofer}
\end{thm}
\begin{proof} The proof will be essentially the same as that of
the Hofer norm once the following continuity lemma for the smooth
case is proved.

\begin{lem} The function $\gamma: Ham(M,\omega) \to \R_+$ is continuous in
the Hamiltonian topology of $Ham(M,\omega)$.
\end{lem}
\begin{proof}
Let $H \mapsto \phi$ and $K \mapsto \psi$. Then the triangle
inequality of $\gamma$ and the inequality $\gamma(\phi) \leq
\|\phi\|$ imply
$$
|\gamma(\phi) - \gamma(\psi)| \leq \gamma(\phi^{-1}\psi) \leq
\|\phi^{-1}\psi\| \leq \|\overline H \# K\|.
$$
In particular, we have \be\label{eq:gammaHK} |\gamma(\phi) -
\gamma(\psi)| \leq \inf_{H\mapsto\phi,\, K \mapsto\psi} \|\overline
H \# K\|. \ee Now let $\phi \in Ham(M,\omega)$ and $\e > 0$ be
given. Recalling the fact that $ev_1$ is an open map (see Corollary
3.17 \cite{oh:hameo1}) we consider the open neighborhood
$ev_1(\CU(\phi_H,\e_1,\e_2))$ of $\phi$ where $\phi_H^1 = \phi$. Now
let $\psi \in ev_1(\CU(\phi_H,\e_1,\e_2))$ i.e., $\psi = \phi_K^1$
for some $\phi_K \in \CU(\phi_H,\e_1,\e_2)$. Then we have
$$
\|\overline H \# K \| \leq \e_1, \quad \overline d(\phi_H,\phi_K)
\leq \e_2
$$
by the definition of $\CU(\phi_H,\e_1,\e_2)$. Therefore if we choose
$\e_1 = \e$ and $\e_2$ is \emph{any finite} number, we have
$|\gamma(\phi) - \gamma(\psi)| < \e$ which proves the continuity of
$\gamma$ in the Hamiltonian topology.
\end{proof}

We omit the rest of the details of the proof referring to the
corresponding proofs of Theorem \ref{extendedHofer}.
\end{proof}

\begin{ques} The following questions seem to be interesting to
study.
\begin{enumerate}
\item Is $\gamma$ (or $\overline\gamma$) continuous in the $C^0$-topology?
\item Does the following identity
\be\label{eq:bargamma=gamma} \overline\gamma|_{Ham(M,\omega)} =
\gamma \ee hold? This is the spectral analog to M\"uller's question,
Question \ref{mueller}.
\end{enumerate}
\end{ques}

\section{Calabi quasi-morphism on $\CP^{ham}(Sympeo(S^2),id)$}
\label{sec:quasimorphism}

In the rest of this section, we will restrict to the case of the
sphere $S^2$ with the standard symplectic form $\omega_{S^2}$ on it.
Omitting the symplectic form $\omega_{S^2}$ from their notations, we
just denote by $\CP^{ham}(Sympeo(S^2),id)$, $Hameo(S^2)$ the groups
of continuous Hamiltonian paths and of Hamiltonian homeomorphisms on
$S^2$ respectively, and so on.

We first state the following proposition which is the path space
version of Theorem 3.1 \cite{entov-pol} by Entov and Polterovich.

\begin{prop}\label{eq:rhoH1}
Consider $S^2$ with the standard symplectic form $\omega_{S^2}$ on
it. Let $H, \, F$ be smooth normalized Hamiltonians satisfying. Then
the spectral invariant $\rho(H;1)$ satisfies \be\label{eq:triangle}
|\rho(\phi_H\phi_F;1) - (\rho(\phi_H;1) + \rho(\phi_F;1))| \leq R
\ee for some constant $R = R(S^2) > 0 $ depending only on
$\omega_{S^2}$ but independent of $H, \, F$. In particular, the map
$$
\rho(\cdot;1): \CP^{ham}(Symp(S^2),id) \to \R
$$
defines a quasi-morphism.
\end{prop}
\begin{proof} The inequality
\be\label{eq:triangle} \rho(\phi_H\phi_F;1) - (\rho(\phi_H;1) +
\rho(\phi_F;1)) \leq 0 \ee is nothing but a special case of the
triangle inequality (\ref{eq:triangle}). The existence of a constant
$R > 0$ such that \be\label{eq:lowerbound} \rho(\phi_H\phi_F;1) -
(\rho(\phi_H;1) + \rho(\phi_F;1)) \geq -R \ee was proved by Entov
and Polterovich (See the proof of Theorem 3.1 \cite{entov-pol} in
the context of the covering space $\widetilde{Ham}(S^2)$ but its
proof equally applies to the context of the path space). Combination
of (\ref{eq:triangle}) and (\ref{eq:lowerbound}) finishes the proof.
\end{proof}
We refer to \cite{gamb-ghys}, \cite{entov-pol} for the general
discussion on the basic properties of the quasi-morphism.

Based on this quasi-morphism $\rho(\cdot; 1)$, Entov and Polterovich
defined a homogeneous quasi-morphism on the universal covering space
$\widetilde{Ham}(S^2,\Omega)$
$$
\widetilde \mu: \widetilde{Ham}(S^2,\Omega) \to \R
$$
by the formula \be\label{eq:tildemu} \widetilde \mu(\widetilde \phi)
= \left(\int_{S^2}\omega_{S^2} \right)\cdot \lim_{i \to \infty}
\frac{\rho(\widetilde \phi^m;1)}{m}. \ee We like to point out that
due to the different conventions used in \cite{entov-pol}, the
negative sign in the equation (17) \cite{entov-pol} does not appear
in our definition.

Obviously this definition of homogeneous quasi-morphism can be
lifted to the level of Hamiltonian paths:

\begin{defn} We define a homogeneous quasimorphism
$$
\mu^{path}:\CP^{ham}(Symp(S^2),id) \to \R
$$
by the same formula \be\label{eq:mulambda} \mu^{path}(\lambda)
=\left( \int_{S^2}\omega_{S^2} \right)\cdot \lim_{i \to \infty}
\frac{\rho(\lambda^m;1)}{m}. \ee
\end{defn}
From the definition above and the hamiltonian-continuity of
$\rho(\cdot;1)$, it follows that $\mu^{path}$ is also
hamiltonian-continuous.

The following two propositions concerning the quasi-morphism
$\mu^{path}$ were essentially proved by Entov and Polterovich
\cite{entov-pol}.

\begin{prop}[Compare with Proposition 3.3 \cite{entov-pol}]
Suppose that $U \subset S^2$ that is displaceable, i.e., there
exists $\phi \in Ham(S^2)$ such that $\phi(\overline U) \cap
\overline U = \emptyset$. Then we have the identity
$$
\mu^{path}(\lambda) = \operatorname{Cal}^{path}(\lambda)
$$
for all $\lambda$ with
$$
\operatorname{supp}\lambda \subset U.
$$
\end{prop}

Entov and Polterovich called this property the \emph{Calabi
property} of a quasi-morphism. We recall that
$\operatorname{Cal}^{path}(\lambda)$ is defined by the integral
\be\label{eq:Calabiinvariant} \operatorname{Cal}^{path}(\lambda) =
\int_0^1 \int_M H(t,x)\, \Omega_\omega \ee when $\lambda = \phi_H$.
Here $\Omega_\omega$ is the Liouville volume form normalized so that
$\int_M \Omega_\omega = 1$.

\begin{prop}[Proposition 3.4 \cite{entov-pol}]
\label{entov-pol} The quasi-morphism $\mu^{path}$ pushes down to a
homogeneous quasi-morphism $\mu:Ham(S^2) \to \R$. Furthermore $\mu$
is continuous on $Ham(S^2)$ with respect to the Hamiltonian
topology.
\end{prop}
\begin{proof} The proof of the first fact verbatim follows from that of
Proposition 3.4 \cite{entov-pol}. The continuity statement
immediately follows from the hamiltonian-continuity of
$\rho(\cdot;1)$ and the definition of the Hamiltonian topology on
$Ham(S^2)$.
\end{proof}

For any given open subset $U \subset S^2$, we denote by
$$
\CP^{ham}(Symp_U(S^2),id)
$$
the set of Hamiltonian paths supported in $U$.

An immediate corollary of these two propositions is the following
homomorphism property of $\mu$ restricted to
$\CP^{ham}(Symp_U(S^2),id)$.

\begin{cor}\label{displacablemu} Suppose that $U$ is an open subset of $S^2$
such that $\overline U$ is displaceable on $S^2$ and let
$$
\lambda_1, \, \lambda_2 \in \CP^{ham}(Symp_U(S^2),id).
$$
Then we have
$$
\mu^{path}(\lambda_1\lambda_2) = \mu^{path}(\lambda_1) +
\mu^{path}(\lambda_2).
$$
\end{cor}

Now we extend all the above discussions to the level of
\emph{continuous} Hamiltonian paths. But these generalization
immediately follow once we know the facts that
\begin{enumerate}
\item $\rho(\cdot;1)$ has been extended to
$\CP^{ham}(Sympeo(M,\omega),id)$ for an arbitrary closed symplectic
manifold, i.e., in particular for $(S^2,\omega_{S^2})$ in section
\ref{sec:spectral}. \item In addition, this extension is
hamiltonian-continuous, i.e, continuous in the Hamiltonian topology.
\end{enumerate}

We summarize the above discussion into the following theorem.

\begin{thm}\label{barmupath} We have an extension of
$\mu^{path}:\CP^{ham}(Symp(S^2),id) \to \R$ to a quasi-morphism
$$
\overline \mu^{path}:\CP^{ham}(Sympeo(S^2),id) \to \R
$$
that satisfies all the analogs to Proposition \ref{entov-pol} and
the Calabi property.
\end{thm}

Now we state the following conjecture, which we strongly believe
would play an essential role in the study of simpleness question of
the area preserving group of $S^2$ (and also of $D^2$). (See Theorem
\ref{simpleness} and \ref{thm:simpleonS2} for some indication.)
Recall from \cite{entov-pol} that the corresponding fact was proved
by Entov and Polterovich for the group $Ham(S^2)$ of \emph{smooth}
Hamiltonian diffeomorphisms on $S^2$.

\begin{conj}\label{barmu} Let $\overline\mu^{path}:\CP^{ham}(Sympeo(S^2),id) \to \R$
be the above extension of the homogeneous Calabi quasi-morphism
given in (\ref{eq:mulambda}). This pushes down to a homogeneous
quasi-morphism $\overline\mu: Hameo(S^2) \to \R$ that satisfies
\be\label{eq:project} \overline\mu^{path} = \overline\mu \circ
\overline{ev}_1. \ee In particular, $\overline{\mu}^{path}(\lambda)$
depends only on the time-one map $\lambda(1)$ of $\lambda$ as long
as $\lambda$ lies in $\CP^{ham}(Sympeo(S^2),id)$.
\end{conj}

An immediate corollary of Conjecture \ref{barmu} and of the Calabi
property of $\overline\mu^{path}$ would be the solution to the
following conjecture

\begin{conj}[Fathi \cite{fathi}] \label{extendedcalabi} The Calabi homomorphism
$\operatorname{Cal}: Ham(D^2,\del D^2) \to \R$ is extended to a
homomorphism
$$
\overline{\operatorname{Cal}}: Hameo(D^2,\del D^2) \to \R
$$
that is continuous in Hamiltonian topology.
\end{conj}

In the next section, we will explain how validity of this conjecture
together with the smoothing theorem \cite{oh:smoothing}, would imply
properness of $Hameo(D^2,\del D^2)$ in $Homeo^\Omega(D^2,\del D^2)$
and hence lead to the proofs of non-simpleness both of
$Homeo^\Omega(D^2,\del D^2)$ and of $Homeo^\Omega(S^2)$.

\section{Discussion: wild area preserving homeomorphisms on $D^2$ and on $S^2$}
\label{sec:wild}

In this section, we will describe an example of a compactly
supported area preserving homeomorphism in $Sympeo(D^2,\partial
D^2)$ that would not be contained in $Hameo(D^2,\partial D^2)$,
\emph{if Conjecture \ref{extendedcalabi} should hold}. Then this
would imply that $Hameo(D^2,\partial D^2)$ is a proper normal
subgroup of $Sympeo(D^2,\partial D^2)$. Combination the above chain
of statements would give rise to non-simpleness of
$Homeo^\Omega(D^2,\del D^2)$, via the following theorem which is a
corollary of the smoothing theorem from \cite{oh:smoothing},
\cite{sikorav}.

\begin{thm}[Theorem I \cite{oh:smoothing}]\label{smoothing} We have
$$
Sympeo(D^2,\partial D^2) = Homeo^\Omega(D^2,\del D^2).
$$
for the standard area form $\Omega$ on $D^2$ regarding it also as
the symplectic form $\omega = \Omega$.
\end{thm}

This being said, we will focus on construction of an example of a
wild area-preserving homeomorphism on $D^2$. For this description,
we will need to consider the conjugate action of rescaling maps of
$D^2$
$$
R_a: D^2(1) \to D^2(a) \subset D^2(1)
$$
for $0 < a < 1$ on $Hameo(D^2,\partial D^2)$, where $D^2(a)$ is the
disc of radius $a$ with its center at the origin. We note that $R_a$
is a conformally symplectic map and so its conjugate action maps a
symplectic map to a symplectic map whenever it is defined.

Furthermore the conjugation by $R_a$ defines a map
$$
\phi \mapsto  R_a^{-1} \circ \phi \circ R_a \,; \,
Hameo(D^2,\partial D^2) \to Hameo(D^2(a),\partial D^2(a)) \subset
Hameo(D^2,\partial D^2)
$$
and the conjugation by $R_a^{-1}$ defines a map
$$
Hameo(D^2(a),\partial D^2(a)) \to Hameo(D^2,\partial D^2).
$$

We have the following important formula for the behavior of Calabi
invariants under the Alexander isotopy.

\begin{lem}\label{Callambdaa} Let $\lambda$ be a given compactly supported
continuous Hamiltonian path on $D^2$ and $\eta > 0$ be a small
constant such that $\mbox{supp }\lambda \subset D^2(1-\eta)$. We
define $\lambda_a: D^2 \to D^2$ defined by
$$
\lambda_a(t,x) =
\begin{cases}a \lambda(t,\frac{x}{a})  \quad &\mbox{for $|x|
\leq a(1-\eta)$} \\
x \quad &\mbox{otherwise}
\end{cases}
$$
for $0 < a \leq 1$. Then $\lambda_a$ is also a continuous
Hamiltonian path on $D^2$ and satisfies \be\label{eq:callamlama}
\overline{\operatorname{Cal}}^{path}(\lambda_a) = a^4
\overline{\operatorname{Cal}}^{path}(\lambda). \ee
\end{lem}
\begin{proof} A straightforward calculation proves that $\lambda_a$
is generated by the (unique) continuous Hamiltonian defined by
$$
\operatorname{Dev}(\lambda_a)(t,x) = \begin{cases} a^2 H \left(t,
\frac{x}{a}\right)
\quad &\mbox{for $|x| \leq a(1-\eta)$} \\
0 \quad &\mbox{otherwise}
\end{cases}
$$
where $H = \operatorname{Dev}(\lambda)$: Obviously the right hand
side function is the hamiltonian-limit of
$\operatorname{Dev}(\lambda_{i,a})$ for a sequence $\lambda_i$ of
smooth hamiltonian approximation of $\lambda$ where  $\lambda_{i,a}$
is defined by the same formula for $\lambda_i$.

From this, we derive the formula \beastar
\overline{\operatorname{Cal}}^{path}(\lambda_a) & = & \int _0^1
\int_{D^2(a(1-\eta))} a^2 H\left(t, \frac{x}{a}\right)
\Omega \wedge dt \\
& = & a^4 \int_0^1 \int_{D^2}H(t,y)\Omega \wedge \, dt = a^4
\overline{\operatorname{Cal}}^{path}(\lambda) \eeastar This proves
(\ref{eq:callamlama}).
\end{proof}

Here comes a construction of an example of wild area preserving
homeomorphisms, which is an enhancement of the one described in
Example 4.2 \cite{oh:hameo1}.

\begin{exm}\label{example}

With the above preparations, we consider the set of dyadic numbers
$\frac{1}{2^k}$ for $k = 0, \cdots$. Let $(r,\theta)$ be polar
coordinates on $D^2$. Then the standard area form is given by
$$
\omega = r\, dr \wedge d\theta.
$$
Consider maps $\phi_k: D^2 \to D^2$ of the form given by
$$
\phi_k = \phi_{\rho_k}: (r,\theta) \to (r, \theta + \rho_k(r))
$$
where $\rho_k: (0,1] \to [0, \infty)$ is a smooth function supported
in $(0,1)$. It follows $\phi_{\rho_k}$ is an area preserving map
generated by an autonomous Hamiltonian given by
$$
F_{\phi_k}(r,\theta) = - \int_1^r s \rho_k(s)\, ds.
$$
Therefore its Calabi invariant becomes \be\label{eq:phikcalabi}
\operatorname{Cal}(\phi_k) = -\int_{D^2}\left(\int_1^r s \rho_k(s)\,
ds\right) r\, dr\,d\theta = 2\pi \int_0^1 r^2 \rho_k(r)\, dt. \ee We
now choose $\rho_k$ in the following way:
\begin{enumerate}
\item $\rho_k$ has support in $\frac{1}{2^k} < r < \frac{1}{2^{k-1}}$
\item For each $k = 1, \cdots$, we have
\be\label{eq:rhokrhok-1} \rho_k(r) = 2^4 \rho_{k-1}(2r) \ee for $r
\in (\frac{1}{2^k},\frac{1}{2^{k-1}})$.
\item $\operatorname{Cal}(\phi_1) = 1$.
\end{enumerate}
Since $\phi_k$'s have disjoint supports by construction, we can
freely compose without concerning about the order of compositions.
It follows that the infinite product
$$
\Pi_{k=0}^\infty \phi_k
$$
is well-defined and defines a continuous map that is smooth except
at the origin at which $\phi_\rho$ is continuous but not
differentiable: This infinite product can also be written as the
homeomorphism having its values given by $\phi_\rho(0) = 0$ and
$$
\phi_\rho(r,\theta)=(r, \theta + \rho(r))
$$
where the smooth function $\rho:(0,1] \to \R$ is  defined by
$$
\rho(r) = \rho_k(r) \quad \mbox{for
$[\frac{1}{2^k},\frac{1}{2^{k-1}}],\, k = 1,\, 2, \cdots$}.
$$
It is easy to check that $\phi_\rho$ is smooth $D^2 \setminus \{0\}$
and is a continuous map, even at $0$, which coincides with the above
infinite product. Obviously the map $\phi_{-\rho}$ is the inverse of
$\phi_\rho$ which shows that it is a homeomorphism. Furthermore we
have
$$
\phi_\rho^*(r\, dr\wedge d\theta) = r\,dr\wedge d\theta \quad
\text{on } \, D^2 \setminus \{0\}
$$
which implies that $\phi_\rho$ is indeed area preserving.
\end{exm}

The following lemma will play an important role in our proof of
Theorem \ref{simpleness}.

\begin{lem}\label{conjugate}
Let $\phi_k$ the diffeomorphisms given in Example \ref{example}. We
have the identity \bea\label{eq:conjugate} R_{\frac{1}{2}}\circ
\phi_{k-1}^{2^4} \circ R_{\frac{1}{2}}^{-1} = \phi_k. \eea In
particular, we have \be\label{eq:Calconjugate}
\operatorname{Cal}(\phi_k) = \operatorname{Cal}(\phi_{k-1}). \ee
\end{lem}
\begin{proof} Using (\ref{eq:rhokrhok-1}), we compute
$$
R_{\frac{1}{2}} \circ \phi_{k-1} \circ
R_{\frac{1}{2}}^{-1}(r,\theta) = (r, \theta + \rho_{k-1}(2r)) =
\left(r,\theta + \frac{1}{2^4}\rho_k(r)\right)
$$
where the second identity follows from (\ref{eq:rhokrhok-1}).
Iterating this identity $2^4$ times, we obtain (\ref{eq:conjugate})
from (\ref{eq:rhokrhok-1}). The equality (\ref{eq:Calconjugate})
follows from this and (\ref{eq:callamlama}).
\end{proof}

An immediate corollary of this lemma and (\ref{eq:rhokrhok-1}) is
the following
\begin{cor}\label{equality}
We have
$$
\operatorname{Cal}(\phi_k) =  1.
$$
for all $k = 1, \cdots$
\end{cor}

Now we are ready to give the proof of the following theorem.

\begin{thm}\label{simpleness} Validity of
Conjecture \ref{extendedcalabi} implies that $\phi_\rho$ cannot be
contained in $Hameo(D^2,\partial D^2)$.
\end{thm}
\begin{proof}
Suppose to the contrary that $\phi_\rho \in Hameo(D^2,\partial
D^2)$.

Then its Calabi invariant has a finite value which we denote
\be\label{eq:Cal=C} \overline{\operatorname{Cal}}(\phi_\rho) = C_1
\ee for some $C_1 \in \R$. We will derive a contradiction out of
this finiteness.

Writing $\phi_\rho = \psi_N \widetilde \psi_N$ where \beastar
\psi_N & = & \Pi_{i=1}^N \phi_i \\
\widetilde \psi_N & = & \Pi_{i=N+1}^\infty \phi_i, \eeastar we
derive \be\label{eq:Cal=C} C_1 =
\overline{\operatorname{Cal}}(\psi_N) +
\overline{\operatorname{Cal}}(\widetilde \psi_N) \ee from the
homomorphism property of $\overline{\operatorname{Cal}}$. Here we
note that $\psi_N$ is smooth and so obviously lies in
$Hameo(D^2,\del D^2)$. Therefore it follows from the group property
of $Hameo(D^2,\del D^2)$ that $\widetilde \psi_N$ lies in
$Hameo(D^2, \del D^2)$ if $\phi_\rho$ does so.

Now we set $N =1$ and derive
$$
\overline{\operatorname{Cal}}(\psi_1) = \operatorname{Cal}(\psi_1) =
\operatorname{Cal}(\phi_1) = 1
$$
from Corollary \ref{equality}, and hence \be\label{eq:CaltildeC-N}
\overline{\operatorname{Cal}}(\widetilde \psi_1) = C_1 - 1. \ee

On the other hand, applying (\ref{eq:conjugate}) iteratively to the
infinite product
$$
\widetilde \psi_1 = \prod_{i = 2}^\infty\phi_i,
$$
we show that $\widetilde \psi_1$ satisfies the identity \be
\widetilde \psi_1(r, \theta) =
\begin{cases}
R_{\frac{1}{2}}\circ \phi_\rho^{2^4} \circ R_{\frac{1}{2}}^{-1}(r,
\theta) &\quad \mbox{for }\, 0 < r
\leq \frac{1}{2}\\
(r,\theta) & \quad  \mbox{for }\, \frac{1}{2} \leq r \leq 1.
\end{cases}
\ee Note the identity
$$
R_{\frac{1}{2}}\circ \phi_\rho^{2^4} \circ R_{\frac{1}{2}}^{-1}(r,
\theta)= (R_{\frac{1}{2}}\circ \phi_\rho \circ
R_{\frac{1}{2}}^{-1})^{2^4}(r, \theta).
$$
This, the homomorphism property of $\overline{\operatorname{Cal}}$
and Lemma \ref{Callambdaa} applied for $a = \frac{1}{2}$ give rise
to \bea \overline{\operatorname{Cal}}(\widetilde \psi_N) & = & 2^4
\overline{\operatorname{Cal}}(R_{\frac{1}{2}}\circ \phi_\rho
\circ R_{\frac{1}{2}}^{-1}) \nonumber \\
& = & 2^4 \left(\frac{1}{2^4}\right)
\overline{\operatorname{Cal}}(\phi_\rho) =
\overline{\operatorname{Cal}}(\phi_\rho) = C_1 \label{eq:tildeCal}
\eea It is manifest that (\ref{eq:CaltildeC-N}) and
(\ref{eq:tildeCal}) contradict to each other. This finishes the
proof.
\end{proof}

\bigskip

Next we prove the following $S^2$ analog to Theorem
\ref{simpleness}. We first prove

\begin{thm}\label{thm:5.7imply5.8} Conjecture \ref{barmu} implies Conjecture \ref{extendedcalabi}.
\end{thm}
\begin{proof} Embedding $D^2 \to S^2$ as the upper hemisphere, we identify
$D^2$ with the upper hemisphere $D^+ \subset S^2$. For any given $
\phi \in Hameo(D^2,\del D^2)$, we extend the map to $S^2$ by putting
the identity map on $S^2 \setminus D^+$ and denote the extended map
on $S^2$ by $\widetilde \phi$. Then we define
$$
\overline{\operatorname{Cal}}(\phi): = \overline{\mu}(\widetilde
\phi).
$$
By Corollary \ref{displacablemu}, the Calabi property of
$\overline\mu^{path}$ and Conjecture \ref{barmu}, it follows that
$\overline{\operatorname{Cal}}$ defines a well-defined homomorphism
which extends the usual Calabi homomorphism $\operatorname{Cal}:
Ham(D^2,\del D^2) \to \R$ to $Hameo(D^2,\del D^2)$. This finishes
the proof.
\end{proof}

Next we prove

\begin{thm}\label{thm:simpleonS2} Validity of Conjecture \ref{barmu} implies
that the group $Hameo(S^2)$ is a proper subgroup of $Homeo(S^2)$,
and hence that $Homeo(S^2)$ is not a simple group.
\end{thm}
\begin{proof} Embedding $D^2 \to S^2$ as the upper hemisphere, we identify
$D^2$ with the upper hemisphere $D^+ \subset S^2$ and extend the
homeomorphism $\phi_\rho$ on $D^2$ to an area preserving
homeomorphism on $S^2$ by the identity on $S^2 \setminus D^+$. We
denote the extension by $\widetilde \phi_\rho$.

We claim $\widetilde \phi_\rho$ is not in $Hameo(S^2)$. We denote by
$\overline{\operatorname{Cal}}_{D^+}$ the extension obtained in
Theorem \ref{thm:5.7imply5.8}. Suppose $\widetilde \phi_\rho$ is in
$Hameo(S^2)$ and so $\overline\mu(\widetilde \phi_\rho)$ has a
finite value. Then by the Calabi property of $\overline\mu$ we have
$$
\overline{\operatorname{Cal}}_{D^+}(\phi_\rho) =
\overline\mu(\widetilde \phi_\rho)
$$
and so $\overline{\operatorname{Cal}}_{D^+}(\phi_\rho)$ must have a
finite value. But this gives rise to a contradiction by the proof of
Theorem \ref{simpleness}. This finishes the proof.
\end{proof}

\emph{In conclusion, Conjecture \ref{barmu} is the one to beat!}

\end{document}